\documentclass[reqno]{amsart}

\usepackage{amsmath,graphicx,amssymb,amscd,amsfonts,amsthm}
   \newtheorem{thm}{Theorem}[section]
     \newtheorem{prop}[thm]{Proposition}
     \newtheorem{lem}[thm]{Lemma}
     \newtheorem{cor}[thm]{Corollary}
     
     \newtheorem*{remark}{Remark}
       \newtheorem*{rems}{Remarks}

\newcommand{\R}{\mathbf R}
\newcommand{\Z}{\mathbf Z}
\newcommand{\C}{\mathbf C}
\newcommand{\HH}{\mathbf H}
\newcommand{\T}{\mathbf T}
\newcommand{\Q}{\mathbf Q}
\newcommand{\UU}{\mathcal U}

\makeatletter
\def\@seccntformat#1{\@ifundefined{#1@cntformat}%
   {\csname the#1\endcsname\quad}
   {\csname #1@cntformat\endcsname}
}
\makeatother

\begin{document}
\title{Projections in rotation algebras and theta functions}
\author{Florin P. Boca}
\address{Department of Mathematics, University of Wales, Singleton Park, Swansea SA2 8PP, UK}
\address{Institute of Mathematics of the Romanian Academy, P.O.Box 1-764, RO-70700, Bucharest, Romania}
\curraddr{School of Mathematics, Cardiff University, Senghennydd Road, Cardiff CF2 4YH, UK}
\email{BocaFP@cardiff.ac.uk}

\dedicatory{Dedicated to Professor Marc A. Rieffel on the occasion of his 60th birthday}

\thanks{$^\ast$ Research supported by an EPSRC Advanced Fellowship and GAR 198/1998}
\thanks{$^{\ast\ast}$ Paper presented at the $C^*$-algebren Meeting,
Mathematische Forschungsinstitut, Oberwolfach, February 1-7, 1998 and the 50th British Mathematical Colloquium, Manchester, April 6-9, 1998}


\begin{abstract}
For each $\alpha \in (0,1)$, $A_\alpha$ denotes the universal $C^*$-algebra
generated by two unitaries $u$ and $v$, which satisfy the commutation relation
$uv=e^{2\pi i\alpha} vu$. We consider the order four automorphism $\sigma$ of
$A_\alpha$ defined by $\sigma (u)=v$, $\sigma (v)=u^{-1}$ and describe a method for
constructing projections in the fixed point algebra $A_\alpha^\sigma$, using
Rieffel's imprimitivity bimodules and Jacobi's theta functions. In the case
$\alpha =q^{-1}$, $q\in \Z$, $q\geq 2$, we give explicit formulae for
such projections and find a lower bound for the norm of the Harper operator $u+u^* +v+v^*$.
\end{abstract}

\maketitle

\medskip

The commutative algebra $C({\mathbf T}^n)$ of continuous functions on
the ordinary $n$-dimensional torus $\T^n =\big\{ (z_1,\ldots,z_n ) ; \vert z_j \vert =1\big\}$
is isomorphic to
 the universal $C^*$-algebra generated by $n$ commuting unitary operators.
A non-commutative $n$-torus $A_\alpha$ is the universal $C^*$-algebra
generated by $n$ unitaries $u_1,\ldots,u_n$ subject to relations
$u_j u_k =e^{2\pi i\alpha_{jk}} u_k u_j$, where $\alpha =(\alpha_{jk} )_{1\leq
j,k\leq n}$ is a skew symmetric matrix with real entries.
In some situations it is convenient to regard $\alpha$
as a real skew bilinear form on
$\Z^n$ defined by $\alpha (e_j,e_k)=\alpha_{jk}$ and $A_\alpha$ as the
twisted group $C^*$-algebra $C^* (\Z^n,\beta)$, where
$\beta :\Z^n \times \Z^n \rightarrow \T$ is a $2$-cocycle
such that $\beta (x,y)\, \overline{\beta (y,x)} =e^{2\pi i \alpha (x,y)}$.

In this paper we only consider the case $n=2$, when $\alpha$ is a real
number and $A_\alpha$ is isomorphic (when $\alpha$ is irrational) to the crossed-product $C^*$-algebra
$C(\T)\propto_{\beta} \Z$, where
$\beta$ is the automorphism of $C(\T)$ defined by
$\beta (\phi) (e^{2\pi it})=\phi (e^{2\pi i(t+\alpha)})$,
$\phi \in C(\T$), $t\in \R$.
The algebra $A_\alpha$, called the rotation algebra by angle $\alpha$,
coincides with the universal $C^*$-algebra generated by two unitaries $u$
and $v$ such that $uv =e^{2\pi i\alpha} vu$. It
is endowed with the canonical faithful tracial state $\tau$ defined by
$\tau \Big( \sum\limits_{m,n} a_{m,n} u^m v^n \Big)  =
a_{0,0}$. The starting point in the study of rotation algebras is
the existence of the Powers-Rieffel projections. They are
projections $e_\alpha$ of trace the fractional part $\{ \alpha\}$ of $\alpha$ in $A_\alpha$, which are of the form
$$e_\alpha =G(u)v +F(u)+\check{G} (u)v^{-1}  ,$$
where $\check{G}(x)=G(-x)$, $x\in \R$ and  $F$ and $G$
 are some smooth functions
on $\R$ (see \cite{R1}). The classical results of Pimsner and Voiculescu
(\cite{PV1},\cite{PV2}) show that $\tau_\ast \big( K_0 (A_\alpha )\big) =
\Z +\Z \alpha =\Z +\Z \tau (e_\alpha )$;
 in particular, for any irrational
numbers $\alpha_1$ and $\alpha_2$, the rotation algebras
$A_{\alpha_1}$ and $A_{\alpha_2}$ are isomorphic if and only if
$\alpha_2 \pm \alpha_1 \in \Z$.
 The Powers-Rieffel projections provide a central class of examples of projective modules in
 A. Connes' noncommutative differential geometry (\cite{Co1,Co2}). They also play
 an important r\^ ole in the recent results (\cite{EE},\cite{Bo}) on the
structure of non-commutative tori.

Each matrix $X=\left( \begin{smallmatrix}  a & b \\ c & d \end{smallmatrix} \right) \in
SL_2 (\Z)$ defines a $^\ast$-automorphism $\sigma_X$ of $A_\alpha $,
acting on the canonical generators $u$ and $v$ of $A_\alpha $
as $\sigma_X (u)=u^a v^b$, $\sigma_X (v)=u^c v^d$. Throughout this
paper, we will denote by $\sigma =\sigma_{\left(\begin{smallmatrix}
0 & 1 \\ -1 & 0 \end{smallmatrix}\right)}
$ the ``Fourier transform" automorphism of $A_\alpha$, acting on its generators
by $\sigma (u)=v$ and $\sigma (v)=u^{-1}$. We also set
$$A_\alpha^\sigma =\big\{ a\in A_\alpha   ; \sigma (a)=a\big\}$$
and
$$H_\alpha^{(n)} =u^n +u^{-n} +v^n +v^{-n}  \in   A_\alpha^\sigma ,
\quad n\in \Z .$$

We notice that for any $\alpha \in \R \setminus
\Q$, the $C^*$-algebra
 $A_\alpha^\sigma$ is generated only by the operators
$H_\alpha =H_\alpha^{(1)}$ and $H_\alpha^{(2)}$ (see the appendix).
The $\mbox{\rm K}$-groups of $A_\alpha^\sigma$ were
computed in the case $\alpha = \frac{p}{q} \in \Q$,
$\gcd(p,q)=1$ in \cite{FW}, where it was shown that
 $\mbox{\rm K}_1 (A_\alpha^\sigma )=0$ and
$\mbox{\rm K}_0 (A_\alpha^\sigma )=\Z^9$ if $q\geq 5$.

The problem of characterizing the spectrum of the self-adjoint operator
$H_\alpha=u+u^* +v+v^*$, or of $H_{\alpha,\lambda} =
u+u^* + \frac{\lambda}{2} \, (v+v^*)$, $\lambda \in\R$, in $A_\alpha$
 is very important in the study of
the quantum Hall effect (\cite{Be}).
 If $E_{\alpha ,\lambda}$ denotes the spectral measure of $H_{\alpha ,\lambda}$,
then $\mu_{\alpha ,\lambda} =\tau E_{\alpha ,\lambda}$
 is a measure with $\operatorname{supp} (\mu_{\alpha ,\lambda}) \subset [
-2-\lambda ,2+\lambda ]$ and since
$\tau$ is faithful, its support coincides with the spectrum of
$H_{\alpha ,\lambda}$. Actually one can gather information on $\operatorname{spec} (H_{\alpha ,\lambda})$
 from the K-theoretical
properties of $A_\alpha$. In this respect the results of Pimsner and
Voiculescu (\cite{PV2}) show that for any irrational
$\alpha$, any $\lambda \in \R$ and any
$t$ which belongs to a gap of $\operatorname{spec} (H_{\alpha ,\lambda})$,
there exists an integer $n$ such that
$$\mu_{\alpha ,\lambda} \big( \operatorname{spec}(H_{\alpha,\lambda})
\cap \chi_{(-\infty ,t]} \big) =\{ n\alpha \} .$$

In the irrational case, the label $n$ coincides with the first Chern character of Connes (\cite{Co1}).
Knowing  more about the
projections of $A_\alpha^\sigma$ and about $\tau \big( {\mathcal P}
 (A_\alpha^\sigma ) \big)$, would presumably provide additional
information on $\operatorname{spec} (H_\alpha )$.

Another important feature of the automorphism $\sigma$ is that it
 implements the Andre-Aubry duality.
 One easily checks that the set
$\sigma_+ (\alpha,\lambda)= \bigcup\limits_\theta \,
\sigma (\alpha ,\lambda, \theta)$, as defined in \cite{AMS},
 coincides with the
spectrum of the operator $H_{\alpha,\lambda} \in A_\alpha$.
Using the fact that
$\sigma$ is an automorphism of $A_\alpha$ and
$$\sigma (H_{\alpha,\lambda} )=\sigma \bigg( u+u^* +
\frac{\lambda}{2} \,  (v+v^* )\bigg) =v+v^* +
\frac{\lambda}{2} \,  (u+u^* )=\frac{\lambda}{2} \,  \bigg(
u+u^* +\frac{4}{\lambda} \cdot \frac{v+v^*}{2} \bigg) =
\frac{\lambda}{2} \,  H_{\alpha,\frac{4}{\lambda}} ,$$
it follows that the operators $H_{\alpha,\lambda}$ and
$\frac{\lambda}{2} \, H_{\alpha,\frac{4}{\lambda}}$ have the
same spectrum. This provides (cf. J. Bellissard) a quick proof of the Andre-Aubry duality
$$\sigma_+ (\alpha,\lambda)=\operatorname{spec} \, (H_{\alpha,\lambda})=
 \frac{\lambda}{2} \, \operatorname{spec} \, (H_{\alpha,\frac{4}{\lambda}} )
=\frac{\lambda}{2} \, \sigma_+ \bigg( \alpha,\frac{4}{\lambda} \bigg).$$

The aim of this paper is to develop a method for constructing projections
in the $C^*$-algebra  $A_\alpha^\sigma$.
 In the first section, we prove that for any $\alpha \in (0,1)$,
$A_\alpha^\sigma$ contains a projection $e_\alpha$
of trace $\alpha$ (the trace is the canonical trace on $A_\alpha$). Although
we make use of Rieffel's formalism for constructing
imprimitivity bimodules between twisted $C^*$-group algebras
 associated with lattices in abelian locally compact
groups (\cite{R2}), the nature of these projections
is different from that of the customary Powers-Rieffel projections (\cite{R1}).
Actually, the projections we construct in $A_\alpha^\sigma$ are related to
the classical Jacobi theta functions (\cite{M},\cite{I}).
An important ingredient in the proof
is the inequality
$$\vartheta (0,it)^2 <\vartheta (0,it)+\vartheta \bigg(
\frac{1}{2}  ,it \bigg),\quad  t>\frac{1}{2}  ,$$
which we prove in
Proposition 1.3, employing the infinite product factorization
of theta functions; we use the customary notation
$$
 \vartheta (z,\tau) =
 \sum\limits_m e^{\pi im^2 \tau +2\pi imz} ,\quad z\in \C,\
\tau \in \HH =\{ w\in \C  ;\operatorname{Im} (w) >0 \}  .$$

Another interesting fact related to this approach is that
 the classical
transformation formula for $\vartheta (x,it)$ (see for example \cite[p.33]{M})
$$\vartheta \bigg( \frac{x}{it}  ,\frac{i}{t} \bigg) =\sqrt{t} \,
e^{\frac{\pi x^2}{t}} \vartheta (x,it)  ,\quad
 x\in \R ,\ t>0,$$
arises as a mere consequence of Rieffel's trace formula (\cite[Thm.3.5]{R2}).

In Section 2 we  prove that if $p,q\in \Z$,
$q\geq 2$,
$0< q\alpha -p <\frac{1}{2q}\, $ and $p$ is a quadratic residue
of $q$, then $A_\alpha^\sigma$ contains
projections of trace $q\alpha -p$. If
$0<p-q\alpha < \frac{1}{2q}$ and $-p$ is a quadratic residue of $q$, then
$A_\alpha^\sigma$ contains a projection of trace $p-q\alpha$.

The third section considers the case $\alpha = q^{-1}$,
$q\in \Z$, $q\geq 2$, in more detail. For example, we prove the following
estimates for the norm of the operator $H_\alpha$, $\alpha=q^{-1}$,
$$\| H_\alpha \| \geq \begin{cases}
\displaystyle 4 e^{-\frac{\pi \alpha}{2}}
\frac{\vartheta \big( \frac{1}{2} ,\frac{i}{2\alpha} \big) \,
\vartheta \big( \frac{i}{2} ,\frac{i}{2\alpha} \big)}{\vartheta \big( 0,
\frac{i}{2\alpha} \big)^2}, &  \mbox{if $q$ is even,}  \vspace{0.2cm} \\  \displaystyle 4 e^{-\frac{\pi \alpha}{2}}
\frac{\vartheta \big( \frac{1}{2} ,\frac{i}{2\alpha} \big) \,
\vartheta \big( \frac{i}{2} ,\frac{i}{2\alpha} \big)-2\,
\vartheta^{odd} \big( \frac{1}{2} ,\frac{i}{2\alpha} \big) \,
\vartheta^{odd} \big( \frac{i}{2} ,\frac{i}{2\alpha} \big)}{
\vartheta \big( 0,
\frac{i}{2\alpha} \big)^2 -2\, \vartheta^{odd} \big( 0,
\frac{i}{2\alpha} \big)^2}, & \mbox{if $q$ is odd,}
\end{cases}$$
where we set
$$\vartheta^{odd}  (z,\tau) =\sum\limits_{m\operatorname{odd}}
e^{\pi im^2 \tau +2\pi imz} =\vartheta (z,\tau) -\vartheta (2z,4\tau )
,\quad z\in \C ,\ \tau \in\HH .$$

These estimates should be compared with $\| H_\alpha\|=4-O( \frac{1}{q_1})$, where
$\alpha =1/(q_1+1/(q_2+1/q_3+\cdots)))$ is the continued fraction associated to $\alpha$, proved
by Helffer and Sj\" ostrand in \cite[Thm.1]{HS}.

We derive
 closed formulae for the projection $e_{\frac{1}{q}}$, expressing
 them as averages of products of operator-valued
$\vartheta_{a,b}$ functions as follows:
\begin{equation*}
e_{\frac{1}{q}} =\begin{cases}
\frac{\displaystyle \sum\limits_{r,s=0}^{q-1} e^{-\frac{\pi irs}{q}}
\vartheta^{(q)}_{\frac{s}{q},\frac{r}{2}} (U_2,\tfrac{iq}{2}) \vartheta^{(q)}_{\frac{r}{q},\frac{s}{2}}
(U_1,\tfrac{iq}{2})}{\displaystyle q\vartheta (U_2^q,\tfrac{iq}{2})\vartheta(U_1^q,\tfrac{iq}{2})},
& \mbox{\rm if $q$ is even,} \vspace{0.3cm} \\
\frac{\displaystyle \sum\limits_{r,s=0}^{2q-1} e^{\frac{\pi irs}{q}+\pi i\big[\frac{r}{q}\big] \big[\frac{s}{q}\big]}
\vartheta^{(2q)}_{\frac{s}{2q},0} (U_2,2iq) \vartheta^{(2q)}_{\frac{r}{2q},0} (U_1,2iq)}{
\displaystyle q\sum\limits_{\varepsilon_1,\varepsilon_2=0}^1 e^{\pi i\varepsilon_1\varepsilon_2}
\vartheta^{(2)}_{\frac{\varepsilon_2}{2},0} (U_2^q,2iq) \vartheta^{(2)}_{\frac{\varepsilon_1}{2},0} (U_1^q,2iq) },
& \mbox{\rm if $q$ is odd,}
\end{cases}
\end{equation*}
where we set for any unitary $U$ and any $\tau\in\HH$, $b\in\R$, $N\in {\mathbf N}^*$, $a\in \Z$, $0\leq a<N$,
\begin{equation*}
\begin{split}
& \vartheta^{(N)}_{\frac{a}{N},b} (U,\tau) =\sum\limits_m e^{\pi i\tau(m+\frac{a}{N})^2 +2\pi i(m+\frac{a}{N})b} U^{Nm+a} ,\\
& \vartheta(U,\tau) =\vartheta^{(1)}_{\frac{0}{1},0} (U,\tau) =\sum\limits_m e^{\pi i\tau m^2} U^m .
\end{split}
\end{equation*}

As a result, if we specialize for example to the case when $q$ is even,
the following identity follows for any integer $k\geq 1$ and $t_1,t_2 \in \R$,
$$
\sum\limits_{m,n=0}^{2k-1} e^{-\frac{2\pi imn}{k}}
\vartheta_{\frac{n}{2k},\frac{m}{2}} ( t_1 , ik)
 \vartheta_{\frac{n}{2k},\frac{m}{2}} (- t_1 , ik)
\vartheta_{\frac{m}{2k},\frac{n}{2}} ( t_2 ,ik)
\vartheta_{\frac{m}{2k},\frac{n}{2}} ( - t_2 ,ik )
=2k \vartheta ( t_1 , ik)^2
\vartheta ( t_2 ,  ik)^2
 , $$
where
$$ \vartheta_{a,b} (z,\tau )=\sum\limits_m
e^{\pi i(m+a)^2 \tau +2\pi i(m+a)(z+b)}  ,\quad z\in \C, \
\tau \in \HH .$$

This should be compared with Riemann's identities for theta functions
(see \cite{M}).

\setcounter{equation}{0}

\section{ Constructing projections of trace $\alpha$ in $A_\alpha^\sigma$}

We start by recalling the framework from \cite{R2}.
Let $M$ be an abelian locally compact group, let
 $\widehat{M}$ denote its topological dual and consider
$G=M\times \widehat{M}$ equipped with the Haar-Plancherel measure.
Let $\beta :G\times G\rightarrow
\mbox{\bf T}$ be the
Heisenberg bicharacter defined by
$$\beta \big( (x_1,y_1), (x_2,y_2)\big) =\langle x_1,y_2\rangle,
\quad x_1,x_2\in M,\ y_1,y_2 \in \widehat{M} ,$$
where $\langle \ ,\ \rangle :M\times \widehat{M}
\rightarrow \T$ is the canonical
pairing of $M$ with $\widehat{M}$.

Set $\overline{\beta} (x,y)=\overline{\beta(x,y)}$
and $\beta^\ast (x,y)=\overline{\beta} (y,x)=
\overline{\beta(y,x)}$, $x\in M\, ,\ y\in \widehat{M}$.
The formula
$$\big( \pi_{(x^\prime,x^{\prime \prime})} f\big) (t)=
\langle t,x^{\prime \prime} \rangle  f(t+x^\prime) ,
\qquad t,x^\prime \in M ,\ x^{\prime \prime} \in \widehat{M}$$
defines a square-integrable projective unitary  representation
$\pi:G\rightarrow \UU(L^2 (M))$
such that
$$\begin{array}{l}
\displaystyle \pi_x \pi_y =\beta (x,y)\ \pi_{x+y}  ,\\  \\
\displaystyle
\pi_x \pi_y =\beta \beta^\ast (x,y)\ \pi_y \pi_x ,\\    \\
\displaystyle
(\pi_x)^\ast =\beta(x,x)\, \pi_{-x} ,\qquad x,y\in G .\end{array}$$

If $D$ is a lattice in $G$, we denote by $\vert  G \slash D \vert$
 its covolume and by $C^* (D,\beta )$ the $C^*$-algebra generated
by $\pi_w$, $w\in D$.
The subgroup
$$D^\bot =\big\{ w\in G  ;  \beta \beta^* (D,w) =1\big\} =
\big\{ w\in G  ;
\beta (x,w)=\beta (w,x), \ \forall \, x\in D\big\} \subset G
$$
is a lattice in $G$. We equip
$D$ with the Haar measure
which assigns mass $\vert G\slash D^\bot \vert ^{-1} =
\vert G\slash D\vert$ to each point
 and $D^\bot$ with the Haar measure
assigning mass one
to each point (see \cite[p.278]{R2}).
The twisted $C^*$-algebra $C^*(D,\beta)$
 acts  on the left on
the space ${\mathcal S}(M)$ of Schwartz functions on $M$ by
\begin{equation}
 af=\int\limits_D a(w)\pi_w f dw =
\vert  G\slash D  \vert \sum\limits_{w\in D} a(w) \pi_w f
  ,\qquad f\in {\mathcal S}(M),\ a\in L^1 (D,\beta ) .
\end{equation}

Replacing as in \cite[p.269]{R2} $\pi_z$ by $\pi_z^*$, we regard
$C^*(D^\bot ,\beta)$ as being generated by $\pi_z^*$ acting on the left
on ${\mathcal S}(M)$. This action commutes with the left action of
$C^*(D,\beta)$ for $\pi_w \pi_z =\pi_z \pi_w$, $w\in D$, $z\in D^\bot$. The
opposite algebra of $C^*(D^\bot,\beta)$ is $C^*(D^\bot,\overline{\beta})$,
which acts on the right on ${\mathcal S}(M)$ by
\begin{equation}
 fb =\int\limits_{D^\bot} b(z) \big( \pi_z ^* f\big)  dz
=\sum\limits_{z\in D^\bot} b(z)\pi_z^\ast f  ,\quad
f\in {\mathcal S}(M),\ b\in L^1 (D^\bot ,\overline{\beta})=L^1 (D^\perp,\beta)^{\operatorname{opp}} .
\end{equation}

Moreover, ${\mathcal S}(M)$
becomes a $C^* (D,\beta)-C^* (D^\bot ,\overline{\beta})$ equivalence
bimodule with respect to the $C^*$-valued inner products
$\langle \ \ \! ,\ \rangle_D:{\mathcal S}(M)\times {\mathcal S}(M)
\rightarrow C^*(D,\beta )$ and
$\langle \ \ \! ,\ \rangle_{D^\bot} :{\mathcal S}(M)\times {\mathcal S}(M)
\rightarrow C^*(D^\bot,\overline{\beta} )$ defined for
any $f,g\in {\mathcal S}(M)$, $w=(w^\prime,w^{\prime \prime})\in D$
and $z=(z^\prime ,z^{\prime \prime} )\in D^\bot$ by
\begin{equation}
\big< f,g\big>_D (w)=\big< f,\pi_w g\big>_{L^2 (M)}
 =
\int\limits_{M} f(s) \,
 \overline{g(s+w^\prime)} \ \overline{\langle s,w^{\prime \prime}
 \rangle}  ds  ,
\end{equation}
\begin{equation}
\big< f,g\big>_{D^\bot} (z)=
\big< \pi_z g,f\big>_{L^2 (M)} =
\int\limits_{M} \overline{f(s)} \, g(s+z^\prime )\,
\langle s,z^{\prime \prime}\rangle
ds .
\end{equation}

If $\tau_D$ and $\tau_{D^\bot}$ are the canonical (normalized)
traces on $C^* (D,\beta)$ and respectively on $C^* (D^\bot ,\overline{\beta})$,
then  according to \cite[Thm.3.5]{R2} we have
\begin{equation}
\tau_D \big( \langle f,g\rangle_D \big) =\vert G\slash D \vert
\, \tau_{D^\bot} \big( \langle g,f\rangle_{D^\bot} \big)  .
\end{equation}

In this paper we are interested in the case $M=\mbox{\bf R}^m \times F$,
 where $F$
is a finite cyclic group.
Since $\langle x,y\rangle =\mbox{\bf e}
\bigg( \sum\limits_{j=1}^m x_j y_j
 \bigg) $ if
$x=(x_1,\ldots,x_m)$, $y=(y_1,\ldots,y_m)\in \R^m$ and
$\big< [n]_q,[m]_q \big> =\mbox{\bf e} \big( \frac{nm}{q} \big)$,
$[n]_q$, $[m]_q \in \mbox{\bf Z}_q$, where $\mbox{\bf e} (t)=\exp (2\pi it)$
for all
$t\in \R$,  we shall identify $M$ with
$\widehat{M}$ such that $\langle x,y\rangle =\langle y,x\rangle$ for all
$x,y\in M=\widehat{M}$. Therefore
$R(w_1,w_2)=(-w_2,w_1)$ defines a group automorphism of $G$.
We notice that $R(D^\bot)=(RD)^\bot$. In the sequel we will assume that $RD=D$.

The Fourier transform
$$({\mathcal F} f)(s)=\int\limits_M f(x) \, \overline{\langle x,s\rangle} dx ,
\qquad f\in {\mathcal S} (M) , \ s\in M$$
extends to a unitary on $L^2 (M)$ such that for all $g\in G$,
\begin{equation}
{\mathcal F} \, \pi_g =\overline{\beta (g,g)}\, \pi_{Rg} {\mathcal F}   .
\end{equation}

Since $\beta (Rw,Rw)=\overline{\beta (w,w)}$, we see that for all
$\xi_1,\xi_2 \in {\mathcal S} (M)$ and $w\in D$,
\begin{equation}
\begin{split}
\langle {\mathcal F} \xi_1 ,{\mathcal F} \xi_2 \rangle_D (Rw) & =
\langle {\mathcal F} \xi_1 ,\pi_{Rw}\, {\mathcal F} \xi_2  \rangle_{L^2 (M)} =
\overline{\beta (w,w)}\, \langle {\mathcal F} \xi_1 ,{\mathcal F} \, \pi_w \xi_2
\rangle_{L^2(M)} \\ &
=\overline{\beta (w,w)} \, \langle \xi_1,\pi_w \xi_2 \rangle_{L^2(M)}
=\overline{\beta (w,w)} \, \langle \xi_1,\xi_2 \rangle_D (w) ,
\end{split}
\end{equation}
which yields for all $\xi_1,\xi_2,\xi_3 \in {\mathcal S}(M)$,
\begin{equation}
\begin{split}
 {\mathcal F} \big( \langle \xi_1 ,\xi_2 \rangle_D \xi_3 \big) &
=\int\limits_D \langle \xi_1 ,\xi_2 \rangle_D (w) \, {\mathcal F} \pi_w \xi_3 \, dw
=\int\limits_D \overline{\beta (w,w)}\, \langle \xi_1 ,\xi_2 \rangle_D (w)\,
\pi_{Rw} {\mathcal F} \xi_3 \, dw \\    &
 =\int\limits_D \langle {\mathcal F} \xi_1 ,{\mathcal F} \xi_2 \rangle_D (Rw) \, \pi_{Rw}
{\mathcal F} \xi_3 \, dw =\langle {\mathcal F} \xi_1 ,{\mathcal F}
 \xi_2 \rangle_D ({\mathcal F} \xi_3)  .
\end{split}
\end{equation}

In a similar way, (1.4) and (1.6) yield for all $\xi_1,\xi_2 \in {\mathcal S}  (M)$
and $z\in D^\bot$,
$$\langle {\mathcal F} \xi_1 ,{\mathcal F} \xi_2 \rangle_{D^\bot} (Rz)  =
\langle \pi_{Rz} {\mathcal F} \xi_2 ,{\mathcal F} \xi_1 \rangle_{L^2(M)} =
\beta (z,z) \langle \pi_z \xi_2 ,\xi_1 \rangle_{L^2 (M)}  =
\beta (z,z) \langle \xi_1,\xi_2\rangle_{D^\bot} (z) ,$$
and further on for all $\xi_1,\xi_2 ,\xi_3 \in {\mathcal S}(M)$,
\begin{equation}
\begin{split}
 {\mathcal F}
 \big( \xi_1 \langle \xi_2 ,\xi_3 \rangle_{D^\bot} \big) &  =
\int\limits_{D^\bot} \langle \xi_2 ,\xi_3 \rangle_{D^\bot} (z)\, \beta (z,z)
{\mathcal F} \pi_{-z} \xi_1  dz
=\int\limits_{D^\bot} \langle {\mathcal F} \xi_2 ,{\mathcal F} \xi_3 \rangle_{D^\bot}
(Rz) \, \pi_{Rz}^* {\mathcal F} \xi_1  dz \\ &
 =({\mathcal F} \xi_1) \langle {\mathcal F} \xi_2 ,{\mathcal F} \xi_3 \rangle_{D^\bot}  .
\end{split}
\end{equation}

According to \cite{C}, (1.8) and (1.9) show that $\Z_4 =\Z\slash
4\Z$ acts on the imprimitivity bimodule
${\mathcal S}(M)$ by $g\xi ={\mathcal F} \xi$, $g\in \Z_4$, inducing
automorphisms $\sigma_D \in \mbox{\rm Aut} \big( C^*(D,\beta)\big)$,
$\sigma_{D^\bot} \in \mbox{\rm Aut} \big( C^* (D^\bot ,\beta )\big)$ such that
for all $\xi_1 ,\xi_2 \in {\mathcal S} (M)$, $a\in C^* (D,\beta )$,
$b\in C^*(D^\bot,\beta)$,
\begin{equation}
 \sigma_D \big( \langle \xi_1 ,\xi_2 \rangle_D \big) =
\langle {\mathcal F} \xi_1 ,{\mathcal F} \xi_2 \rangle_D  ,
\end{equation}
\begin{equation}
 {\mathcal F} (a\xi_1)=\sigma_D (a)\, ({\mathcal F} \xi_1)  ,
\end{equation}
\begin{equation}
 \sigma_{D^\bot} \big( \langle \xi_1 ,\xi_2 \rangle_{D^\bot} \big) =
\langle {\mathcal F} \xi_1 ,{\mathcal F} \xi_2 \rangle_{D^\bot}  ,
\end{equation}
\begin{equation}
 {\mathcal F} (\xi_1 b)=({\mathcal F} \xi_1 )\, \sigma_{D^\bot} (b)  .
\end{equation}

To find $\sigma_D$, remark that (1.10) and (1.7) yield
$$\sigma_D \big( \langle \xi_1 ,\xi_2\rangle_D \big) =
\int\limits_D \langle {\mathcal F} \xi_1 ,{\mathcal F} \xi_2 \rangle_D (Rw)\, \pi_{Rw}  dw
=\int\limits_D \langle \xi_1 ,\xi_2 \rangle_D (w) \,
\overline{\beta (w,w)} \, \pi_{Rw}  dw  .$$

On the other hand
$$\sigma_D \big( \langle \xi_1 ,\xi_2 \rangle_D \big) =
\int\limits_D \langle \xi_1 ,\xi_2 \rangle_D (w)\, \sigma_D (\pi_w )  dw,
$$
hence for all $w\in D$,
\begin{equation}
\sigma_D (\pi_w )=\overline{\beta (w,w)}\, \pi_{Rw}  .
\end{equation}

A similar computation shows that for all $z\in D^\bot$,
\begin{equation}
\sigma_{D^\bot} (\pi_z )=\overline{\beta (z,z)} \, \pi_{Rz}  .
\end{equation}

We notice that $\sigma_D^2 (\pi_w)=\pi_{-w}$, $w\in D$,
$\sigma_{D^\bot}^2 (\pi_z )=\pi_{-z}$, $z\in D^\bot$,
$\sigma_D^4 =id_{C^*(D,\beta )}$ and
$\sigma_{D^\bot}^4 =id_{C^*(D^\bot ,\beta)}$.

\begin{prop}
Let $\alpha \in (0,1)$ and
$D=\Z\varepsilon_1 +\Z\varepsilon_2$ be a
lattice in {\em $G=\R^2$} such that $R\varepsilon_1 =\varepsilon_2$ {\em (}so
 $R\varepsilon_2 =-\varepsilon_1${\em )}  and
 $\beta(\varepsilon_j,
\varepsilon_j)=1$, $j=1,2$, $\beta \beta^* (\varepsilon_1,\varepsilon_2)=
e^{2\pi i\alpha}$. Set $U_j =\pi_{\varepsilon_j}$, $j=1,2$ {\em (}hence
$\sigma =\sigma_D$ is an automorphism of $A_\alpha =C^*(D,\beta)$ such that
$\sigma (U_1)=U_2$, $\sigma(U_2)=U_1^{-1}${\em )}.
Assume that there exists $ \lambda \in \{ \pm 1,\pm i \}$ and {\em
$f\in {\mathcal S} (\R)$} such that ${\mathcal F} f=\lambda f$ and
the element
$a=\big< f,f \big> _{D^\bot}$ is invertible in
$C^*(D^\bot,\beta)$.
Then $p=\langle fa^{-1/2},fa^{-1/2} \rangle_D$ is a projection in
$A_\alpha$ such that $\tau_D (p)=\vert G \slash D \vert =\alpha$
 and $\sigma (p)=p$.
\end{prop}

\begin{proof}
The first part follows from \cite{R1}, \cite{R2}, so we only have to prove that
$\sigma (p)=p$. Since ${\mathcal F} f=\lambda f$ and $R(D^\bot)=D^\bot$,
 (1.12) yields $\sigma_{D^\bot} (a)=a$,
 hence $
\sigma_{D^\bot} \big( a^{-1/2} \big) =a^{-1/2}$.
By (1.13)
${\mathcal F} \big( fa^{-1/2} \big)
=\lambda f a^{-1/2}$ and applying (1.10) we get $
\sigma_D (p)=p$.
\end{proof}

Next, we choose $D=\Z \varepsilon_1 +\Z\varepsilon_2$, with
$\varepsilon_1 =\big( \sqrt{\alpha},0 \big)$, $\varepsilon_2 =
\big( 0,\sqrt{\alpha} \, \big)$.
The lattice $D$ has covolume $\vert G\slash D \vert
=\alpha$ in $G=\R^2$ and $RD=D$. Since
$\beta \beta^* (\varepsilon_1,\varepsilon_2)=e^{2\pi i\alpha}$,
the unitaries $U_j =\pi_{\varepsilon_j}$  satisfy
 $U_1 U_2 =e^{2\pi i\alpha} \,  U_2 U_1$, hence
$C^*(D,\beta)$ is canonically isomorphic to $A_\alpha$. The automorphism
$\sigma =\sigma_D$ acts on $A_\alpha =C^*(D,\beta)$ as
$\sigma \big( \pi_{m_1 \varepsilon_1 +m_2 \varepsilon_2} \big) =
e^{2\pi im_1 m_2 \alpha}
\pi_{-m_2 \varepsilon_1 +m_1 \varepsilon_2}$, $m_1,m_2 \in\Z$,
so $\sigma(U_1)=U_2$ and $\sigma (U_2)=U_1^{-1}$.
Notice also that $\sigma^2
\big( \pi_{m_1 \varepsilon_1 +m_2 \varepsilon_2} \big) =
\pi_{-m_1 \varepsilon_1 -m_2 \varepsilon_2} $.
The orthogonal lattice of $D$ with respect to $\beta \beta^*$ is
$D^\bot =\Z \delta_1 +\Z \delta_2$, with $\delta_1 =\big( 0 ,
 \frac{1}{\sqrt{\alpha}} \big)
$, $\delta_2 =\big( \frac{1}{\sqrt{\alpha}}  ,0
\big)$. It has covolume $\alpha^{-1}$
in $\mbox{\bf R}^2$ and $R(D^\bot)=D^\bot$. If we denote $V_1 =\pi^*_{\delta_1}=
\pi_{-\delta_1}$
and $V_2=\pi^*_{\delta_2} =\pi_{-\delta_2}$, then
$\sigma_{D^\bot} (V_1)=V_2^{-1}$ and $\sigma_{D^\bot}(V_2)=V_1$.

Consider also the Schwartz
function $f(s)=e^{-\pi s^2} $, $s\in \R$. Since ${\mathcal F}
 f=f$, Proposition 1.1
 shows that if $a=\langle f,f\rangle_{D^\bot}$ were invertible, then
$  \left< f \langle f, f\rangle_{D^\bot}^{-1/2} ,
f \langle f, f\rangle_{D^\bot}^{-1/2}  \right> _D$ would be a projection of
trace $\alpha$ in $A_\alpha^\sigma$.

The  following formula
 (\cite[p.5]{I}) will be used repeatedly throughout the paper
\begin{equation}
\int\limits_{\mbox{\bf R}} e^{-2\pi s^2 +2\pi as} ds =\frac{1}{\sqrt{2}}\,
e^{ \frac{\pi a^2}{2} }  ,\quad   a\in \C  .
\end{equation}

Note also that for all $t>0$,
$$\int\limits_{\R} e^{-\pi ts^2} ds =\frac{1}{\sqrt{t}}  . $$

We make use of (1.4) and (1.16) to get
\begin{equation}
\begin{array}{l}
\displaystyle
\langle f,f\rangle_{D^\bot} (m_1\delta_1 +m_2 \delta_2)
 =\langle f,f\rangle_{D^\bot} \bigg( \frac{m_2}{\sqrt{\alpha}}\, ,
\frac{m_1}{\sqrt{\alpha}}
\bigg) =\int\limits_{\mbox{\bf R}} \overline{f(s)}\ f\bigg( s+
\frac{m_2}{\sqrt{\alpha}} \bigg)
e^{\frac{2\pi ism_1}{\sqrt{\alpha}}} ds \\  \\
\displaystyle
 \qquad \qquad
 =\int\limits_{\mbox{\bf R}} e^{-2\pi s^2 -\frac{2\pi m_2 s}{\sqrt{\alpha}} +
\frac{2\pi im_1 s}{\sqrt{\alpha}} -\frac{\pi m_2^2}{\alpha} }  ds
 =
 \frac{1}{\sqrt{2}} e^{ -\frac{\pi (m_1^2 +m_2^2)}{2\alpha} -
\frac{\pi im_1 m_2}{\alpha} } ,
\end{array}
\end{equation}
therefore
\begin{equation}
\begin{split}
 a  =\langle f,f\rangle_{D^\bot} & =\sum\limits_{z\in D^\bot}
\langle f,f\rangle_{D^\bot} (z)\pi_z^*
 =\sum\limits_{z\in D^\bot}
\langle f,f\rangle_{D^\bot} (-z)\beta (z,z)\pi_z \\  &
 =\sum\limits_{m_1,m_2} \langle f,f\rangle_{D^\bot}
(-m_1 \delta_1 -m_2 \delta_2)
e^{\frac{2\pi im_1 m_2}{\alpha}}
\pi_{m_1 \delta_1 +m_2 \delta_2} \\ &
 =\frac{1}{\sqrt{2}}
\sum\limits_{m_1,m_2} e^{-\frac{\pi (m_1^2 +m_2^2)}{2\alpha} +
\frac{\pi im_1 m_2}{\alpha} }  V_1^{m_1} V_2^{m_2} \in C^*(D^\perp,\overline{\beta}).
\end{split}
\end{equation}

A direct computation based on (1.18) yields
$ 2\tau_{D^\bot} (\langle f,f\rangle_{D^\bot}^2) =
 \vartheta \big(
0,\frac{i}{\alpha} \big)^2$.
A computation analogue to (1.17)  yields
$\sqrt{2} \langle f,f\rangle_D (m_1 \varepsilon_1 +m_2 \varepsilon_2) =
e^{-\frac{\pi (m_1^2 +m_2^2)\alpha}{2} +\pi im_1 m_2 \alpha }$,
therefore
\begin{equation}
\sqrt{2}  \langle f,f\rangle_D  =
\vert G\slash D \vert  \, \sum\limits_{w\in D} \langle f,f\rangle_D (w)
 \pi_w
 =\alpha \sum\limits_{m_1,m_2}
e^{ -\frac{\pi (m_1^2 +m_2^2)\alpha}{2} +\pi im_1 m_2 \alpha }
U_2^{m_2} U_1^{m_1}  ,
\end{equation}
which yields further
$2\tau_D \big( \langle f,f\rangle_D^2 \big) =\alpha^2 \vartheta
 (0,i\alpha )^2$. On the other hand, (1.5) yields for all
$\phi ,\psi \in {\mathcal S} (\R)$,
\begin{equation}
\tau_D \big( \langle \phi ,\psi \rangle_D \big) =\vert G\slash D \vert
\tau_{D^\bot}
\big( \langle \psi ,\phi \rangle_{D^\bot} \big) =\alpha  \tau_{D^\bot}
\big( \langle \psi ,\phi \rangle_{D^\bot} \big)  ,
\end{equation}
therefore for all $f_1,f_2,f_3,f_4 \in {\mathcal S}(\R)$
\begin{equation}
 \begin{split}
\tau_D \big( \langle f_1,f_2\rangle_D
\langle f_3,f_4 \rangle_D\big)
&  =\tau_D \big( \big< \langle f_1,f_2\rangle_D f_3,f_4 \big>_D \big)
 =\alpha
\tau_{D^\bot} \big( \big< f_4,\langle f_1,f_2\rangle_D f_3\big> _{D^\bot}
 \big) \\  &
=\alpha
 \tau_{D^\bot} \big( \big< f_4, f_1\rangle_{D^\bot}  \langle f_2 ,f_3
\rangle_{D^\bot} \big) .
\end{split}
\end{equation}

Taking $f_j =f$ in the previous equality we get
\begin{equation}
 \vartheta (0,i\alpha)=\frac{1}{\sqrt{\alpha}}\, \vartheta
 \bigg( 0 ,\frac{i}{\alpha} \bigg)  .
\end{equation}

This is one of the modularity conditions satisfied by theta functions.
Its appearance is not really surprising,
for the Poisson summation formula plays
an important r\^ ole in the proof of (1.5) (and implicitly of (1.20)).
 Actually we can do better by taking
$ f_a (s)=e^{ -\pi (s+a)^2 }$, $a\in \R$. A computation similar to
the previous ones gives
\begin{equation*}
\begin{split}
 \tau_D \big( 2 \langle f_a,f\rangle_D  \langle f,f_a \rangle_D \big)
& =\alpha^2  e^{-\pi a^2} \vartheta (-ia\sqrt{\alpha}  ,i\alpha)
\vartheta (0,i\alpha )  , \\
\tau_{D^\bot} \big( 2
 \langle f_a,f_a \rangle_{D^\bot}  \langle f,f\rangle_{D^\bot}
\big) & =\vartheta \bigg( \frac{a}{\sqrt{\alpha}} ,\frac{i}{\alpha} \bigg)
\vartheta \bigg( 0 ,\frac{i}{\alpha} \bigg)  ,
\end{split}
\end{equation*}
and using (1.21) and (1.22),
$$\alpha  e^{-\pi a^2} \vartheta (-ia\sqrt{\alpha}  ,i\alpha )
\vartheta (0,i\alpha) =\vartheta \bigg( \frac{a}{\sqrt{\alpha}}  ,
\frac{i}{\alpha} \bigg)
\vartheta \bigg( 0 ,\frac{i}{\alpha} \bigg) =\frac{1}{\sqrt{\alpha}}
\vartheta \bigg( \frac{a}{\sqrt{\alpha}} ,\frac{i}{\alpha} \bigg)
 \vartheta (0,i\alpha),$$
hence for all $a\in \R$ and $\alpha >0$,
$$\vartheta (-ia\sqrt{\alpha} ,i\alpha)=\frac{1}{\sqrt{\alpha}}\,
e^{\pi a^2}
 \vartheta \bigg( \frac{a}{\sqrt{\alpha}}  ,\frac{i}{\alpha} \bigg).$$

Taking $x=\frac{a}{\sqrt{\alpha}}$, $t=
\frac{1}{\alpha}$, we recover the basic transformation
 formula for $\vartheta$
\begin{equation}
\vartheta \bigg( \frac{x}{it}  ,\frac{i}{t} \bigg) =\sqrt{t}\,
e^{\frac{\pi x^2}{t}} \vartheta (x,it) ,\quad
  x\in \R,\   t>0 .
\end{equation}

\begin{lem}
The operator $X= \sum\limits_m
 e^{-\pi m^2 \alpha_0 }V_1^m$ is positive and
 invertible
for all $\alpha_0 >0$.
\end{lem}

\begin{proof}
 As
 $X \in C^*(V_1)= C(\T)$, we have
$X(\lambda )=\sum\limits_m e^{-\pi m^2 \alpha_0} \lambda^m$,
 $\lambda \in \T$, hence the spectrum of $X$ coincides with
 $\vartheta \big( [0,1],i\alpha_0 \big)$, where
$$ \vartheta (z,\tau) =\sum\limits_m e^{\pi im^2 \tau +2\pi imz} ,
 \quad z\in \C ,\
\tau \in \HH  =\{ \zeta \in \C  ; \operatorname{Im} \zeta >0\}$$
denotes the customary
theta function. The operator $X$ is self-adjoint for
 $$ X(e^{2\pi it} ) =1+
2 \sum\limits_{m\geq 1} e^{-\pi m^2 \alpha_0}\cos (2\pi mt) \in \R
, \quad  t\in \R .$$

On the other hand $X(1)=1+2
\sum\limits_{m\geq 1} e^{-\pi m^2 \alpha_0} >0$, so we only have to show that
$0\notin \vartheta \big( [0,1],i\alpha_0 \big)$. This is true for
$\frac{1+i\alpha_0}{2}$ is the only zero of the function
$\theta_{\alpha_0} (z)=\vartheta (z,i\alpha_0)$ in the fundamental domain
\begin{equation*}
\{ z\in \C ;
0\leq \operatorname{Re} z< 1,\ 0\leq \operatorname{Im} z< \alpha_0 \}  .\qedhere
\end{equation*}
\end{proof}

Let $\alpha \in (0,\infty )$. We set
\begin{equation*}
\begin{split}
 \alpha_m & =e^{-\frac{\pi m^2}{2\alpha} } ,\\
 \beta_m & =\vartheta \bigg( -\frac{m}{2\alpha}  ,
\frac{i}{2\alpha} \bigg) ,\quad m\in \Z ,\\
 c(t) & =\min\limits_{x\in \R} \, \vartheta ( x,it) >0 ,\quad  t>0 ,\\
 C(t) & = \max \limits_{x\in \R}
\, \vartheta ( x,it) =\vartheta ( 0,it)
=\sum\limits_n e^{-\pi n^2 t}  ,\quad t>0 ,\\
 \Phi_m (t) & =\vartheta \bigg( -\frac{t}{\sqrt{\alpha}} +\frac{m}{2\alpha}  ,
\frac{i}{2\alpha} \bigg) ,\quad m\in \Z .
\end{split}
\end{equation*}

We have $\vert \Phi_m (t)\vert
\leq  \sum\limits_{n} e^{-\frac{\pi n^2}{2\alpha}}
 =\vartheta \big( 0 ,\frac{i}{2\alpha} \big) =
C\big( \frac{1}{2\alpha} \big)$ for all $m\in \Z$ and $\alpha >0$,
 hence
the multiplication operator $D_{m} =M_{\Phi_m}$
 is bounded on $L^2 (\R)$ and
$\| D_{m} \| =\| \Phi_m \|_\infty
 \leq C\big( \frac{1}{2\alpha} \big)$.
 On the other hand $\Phi_0 (t)=\vartheta \big( -
\frac{t}{\sqrt{\alpha}}
 ,\frac{i}{2\alpha} \big) \in \R$ and $\Phi_0
 (t)\geq c\big( \frac{1}{2\alpha} \big) >0$ for all $t\in \R$,
 hence $D_{0}$ is invertible and for all $m\in \Z$,
\begin{equation}
\| D_{0}^{-1} D_{m} \| =\big\| M_{\frac{\Phi_m}{\Phi_0}} \big\| =
\sup\limits_{t\in \R}
 \, \bigg|
\frac{\Phi_m (t)}{\Phi_0 (t)} \bigg| \leq \frac{C\big(
 \frac{1}{2\alpha}
\big)}{c\big( \frac{1}{2\alpha} \big)}  .
\end{equation}

We need more information  on the behaviour of the function $
 \vartheta (\cdot ,it)$ on $\R$ and prove the following

\begin{prop}
{\em (i)}
For any $t>0$
$$c(t)=\vartheta \bigg( \frac{1}{2}  , it\bigg).$$

{\em (ii)} For any $t>0.531$
$$\frac {C(t)\big( C(t)-1\big)}{c(t)} <1 \, .$$
\end{prop}

\begin{proof}
The proof relies on the infinite product expansion for theta functions
 (\cite[Prop.14.1]{M}),
which says that for all $z\in \C$ and $\tau \in \HH$,
\begin{equation}
\vartheta (z,\tau) =\prod\limits_{m\geq 1} \big( 1-e^{2\pi im\tau}\big)
\prod\limits_{m\geq 0} \big( 1+ e^{ (2m+1)\pi i\tau -2\pi iz} \big)
\big( 1+ e^{(2m+1)\pi i\tau +2\pi iz} \big)  .
\end{equation}

Actually we will only use the easier fact that
\begin{equation}
\vartheta (z,\tau) =k_\tau
\prod\limits_{m\geq 0} \big( 1+ e^{(2m+1)\pi i\tau -2\pi iz}\big)
\big( 1+ e^{(2m+1)\pi i\tau +2\pi iz} \big) ,
\end{equation}
for some constant $k_\tau \neq 0$ which does not depend on $z$.
To prove (i), we remark that for any $r_m \geq 0$ and any
$\rho \in \T$ we have
\begin{equation}
(1+r_m \rho )\, (1+r_m \overline{\rho} ) \geq (1-r_m )^2  .
\end{equation}

Taking $r=e^{-\pi t} \in (0,1)$, $t>0$, $r_m  =r^{2m+1}$ and $\rho =e^{2\pi ix}$,
 $x\in \R$, we obtain from (1.26) and (1.27),
$$\frac{\vartheta (x,it)}{\vartheta \big( \frac{1}{2} ,it\big)}
=\prod\limits_{m\geq 0} \frac{(1+r_m \rho)(1+r_m \overline{\rho})}{
(1-r_m)^2}  \geq  1 .$$

Since $\vartheta (x,it)>0$, $x\in \R$, (i) follows.

To prove (ii), we fix
 $t>0$ and denote $ P=\sqrt{\frac{C(t)}{c(t)}}=
\sqrt{\frac{\vartheta (0,it)}{\vartheta (
\frac{1}{2}  ,it)}} . $ Equality (1.26) yields
$$\ln P =\sum\limits_{m\geq 0} \big(  \ln  (1+r^{2m+1})-\ln
(1-r^{2m+1}) \big)  .$$

The mean value theorem yields for any $\varepsilon \in (0,1)$,
$$\frac{2\varepsilon}{1+\varepsilon } <\ln (1+\varepsilon)-\ln (1-\varepsilon )
<\frac{2\varepsilon}{1-\varepsilon}   ,$$
hence
$$\ln P <2\sum\limits_{m\geq 0} \frac{r^{2m+1}}{1-r^{2m+1}} <\frac{2r}{1-r}
\sum\limits_{m\geq 0} r^{2m} =\frac{2r}{(1-r)(1-r^2)}$$
and therefore
\begin{equation}
\frac{C(t)}{c(t)} =P^2   <   h(t)=\exp
\left( \frac{4e^{-\pi t}}{(1-e^{-\pi t})(1-e^{-2\pi t})} \right)
= \exp \left( \frac{4r}{(1-r)(1-r^2)} \right)  .
\end{equation}
Using also
$$C(t)-1 =\sum\limits_{n\neq 0} e^{-\pi n^2 t}
=2\sum\limits_{n\geq 1} r^{n^2} <
2\sum\limits_{n\geq 0} r^{3n+1} =\frac{2r}{1-r^3}  ,$$
we get
$$\frac{C(t)\big( C(t)-1\big)}{c(t)} < g(r) =\frac{2r}{1-r^3} \cdot
\exp \left( \frac{4r}{(1-r)(1-r^2)} \right)  .$$

The derivative of
$$ \phi(r)=\ln g(r)=
 \frac{4r}{(1-r)(1-r^2)} +\ln 2 +\ln \bigg(
\frac{r}{1-r^3} \bigg) \, ,\quad r\in (0,1)  ,$$
is
$$\phi^\prime (r) =\frac{4}{(1-r)(1-r^2)} +
\frac{4r(3r+1)}{(1-r)(1-r^2)^2}
+\frac{2r^3+1}{r(1-r^3)} >0  ,\quad \forall \,
r\in (0,1),$$
therefore $\phi$ is monotonically increasing on $(0,1)$.
Moreover, the equation $\phi (r)=0$ has a unique solution
$r_0 \approx 0.188948$ in $(0,1)$. If we set
$\psi (t)=\phi ( e^{-\pi t}),$ the function $\psi$ is then monotonically decreasing
on $(0,\infty)$ and has a unique root $t_0=-\frac{1}{\pi} \log r_0 \approx 0.53094$, therefore
\begin{equation*}
C(t) \big( C(t)-1\big) < c(t),\quad \mbox{\rm for all $t>t_0 =0.530394$.}\qedhere
\end{equation*}
\end{proof}

\begin{prop}
Let $\alpha \in ( 0,0.942]$. Then
$a_0 =\sum\limits_m \alpha_m D_m V_2^m$
 is  a bounded invertible operator in $C^* (D^\bot ,\overline{\beta} )=A_{\frac{1}{\alpha}}$
 and $a_0=\sqrt{2}  \langle f,f\rangle_{D^\bot}$.
\end{prop}

\begin{proof}
Since $\frac{1}{2\alpha} \geq 1>0.531$, (1.24) and the
previous proposition yield
$$\sum\limits_{m\neq 0} \alpha_m \| D_0^{-1} D_m V_2^m \|
\leq \frac{C\big( \frac{1}{2\alpha}\big)}{c\big( \frac{1}{2\alpha} \big)}
 \sum\limits_{m\neq 0} \alpha_m
 \leq  \frac{C\big( \frac{1}{2\alpha} \big) \big( C\big( \frac{1}{2\alpha} \big)-1
\big)}{c\big( \frac{1}{2\alpha} \big)} <1 .$$

This shows that
$I+\sum\limits_{m\neq 0} \alpha_m D_0^{-1} D_m V_2^m$
defines a bounded invertible
operator, hence so is
$ a_0 = \sum\limits_m \alpha_m D_m V_2^m$.
 The operators $a_0$ and $\sqrt{2}  \langle f,f\rangle_{D^\bot} \in
 {\mathcal B} (L^2 (\R))$
coincide. To see this, notice that for all $m_1 ,m_2 \in \Z$,
$$\big( V_1^{m_1} V_2^{m_2} \phi \big) (s)=e^{-\frac{2\pi im_1 s}{\sqrt{\alpha}}}
\phi \bigg( s-\frac{m_2}{\sqrt{\alpha}} \bigg)  ,$$
hence we may use (1.18) to obtain for all $\phi \in L^2 (\R)$,
$s\in \R$,
\begin{equation*}
\begin{split}
\sqrt{2} \big( \langle f,f\rangle_{D^\bot} (\phi )\big) (s) & =
\sum\limits_{m_1,m_2} e^{-\frac{\pi (m_1^2 +m_2^2)}{2\alpha} +\frac{\pi im_1 m_2}{\alpha}
-\frac{2\pi im_1 s}{\sqrt{\alpha}}} \phi \bigg( s-\frac{m_2}{\sqrt{\alpha}} \bigg)
\\
& =\left( \sum\limits_{m_2} \alpha_{m_2} D_{m_2} V_2^{m_2} \phi \right) (s) =
(a_0 \phi )(s).\qedhere
\end{split}
\end{equation*}
\end{proof}

\begin{cor}
For any $\alpha \in (0,1)$, the rotation
algebra $A_\alpha$
contains a projection $e=e_\alpha$ of trace $\alpha$ such that
$\sigma (e)=e$.
\end{cor}

\begin{rems}
{\em (i)
 For all $m\in \Z$ we have $D_m V_2 =V_2 D_{m+2}$.}

{\em (ii)
If $\alpha$ is irrational, the rotation algebra $A_{\frac{1}{\alpha}}
=C^*(D^\bot ,\overline{\beta})$ is simple, thus is isomorphic to the $C^*$-algebra {\em
$C^*(W_1,W_2) \subset {\mathcal B} (\ell^2 (\Z))$} generated by the unitaries
$W_1 \xi_k =\xi_{k+1}$,
$W_2 \xi_k =e^{\frac{2\pi ik}{\alpha}} \xi_k$, where $\xi_k$ is an
orthonormal basis of $\ell^2 (\Z)$.
 If $\rho =e^{\frac{2\pi i}{\alpha}}$, then
the matrix coefficients of $a=\langle f,f\rangle_{D^\bot}$
 in this representation of $A_{\frac{1}{\alpha}}$ are
\begin{equation*}
\begin{split}
\langle a\xi_k ,\xi_l \rangle
&  =\frac{1}{\sqrt{2}}  \sum\limits_{m_1,m_2} \rho^{\frac{m_1 m_2}{2}} \alpha_{m_1}
\alpha_{m_2} \langle W_1^{m_1} W_2^{m_2} \xi_k ,\xi_l \rangle \\ &
 =\frac{1}{\sqrt{2}}  \sum\limits_{m_1,m_2} \rho^{\frac{m_1 m_2}{2} +m_2 k}
\alpha_{m_1} \alpha_{m_2} \langle \xi_{k+m_1},\xi_l \rangle  \\
&
 =\frac{1}{\sqrt{2}}\, \alpha_{l-k} \sum\limits_{m_2} \rho^{\frac{(k+l)m_2}{2}}
\alpha_{m_2}
 =\frac{1}{\sqrt{2}} \alpha_{l-k}\ \vartheta \bigg( \frac{k+l}{2\alpha} ,
\frac{i}{2\alpha} \bigg) \\  &
 =\frac{1}{\sqrt{2}}\, \alpha_{l-k} \beta_{l+k},\qquad  k,l\in \Z
,\end{split}
\end{equation*}
therefore $a\sqrt{2} =\sum\limits_m \alpha_m D_m U^m $, where
$U\xi_k =\xi_{k+1}$, $D_m \xi_k =\beta_{2k-m} \xi_k$,
$m,k\in \Z$. The diagonal
 operators are bounded and invertible because $0<c\leq \beta_n \leq C$ for
all $n\in \Z$.}
\end{rems}

\setcounter{equation}{0}

\section{Existence of projections of trace $q\alpha  -p$ and $p-q\alpha $}

Let $q\in {\mathbf N}^*$, $q\geq 2$ and $p\in \Z$ such that
$ 0<\gamma =\alpha - \frac{p}{q} \leq \frac{1}{2}$ and there exists
$p_0 \in \Z$ such that $p=p_0^2 \mod{q}$. We choose
 $M=\R \times \Z_q$ and $D=\Z \varepsilon_1 +
\Z \varepsilon_2
\subset G =M\times \widehat{M}$, with
$$
 \varepsilon_1 =\big( \sqrt{\gamma},[p_0]_q ,0,[0]_q \big)
\quad \mbox{\rm and} \quad
\varepsilon_2 =\big( 0,[0]_q ,\sqrt{\gamma} ,[p_0]_q \big)  .
$$

Then $D$ is a lattice in $G$ and $ [0,\sqrt{\gamma} \, )\times \Z_q \times
[0,\sqrt{\gamma} )\times \Z_q$ a fundamental domain for $ G\slash D$.
Since $\Z_q$ is endowed with the Haar-Plancherel measure (which assigns
mass $q^{-1/2}$ to each point from $\Z_q$), we get $\vert G\slash D \vert =
q\gamma =q\alpha -p$. An easy computation gives $ D^\bot =\Z \delta_1
+\Z \delta_2$ with
$$\delta_1 =\bigg( 0,[0]_q ,\frac{1}{q\sqrt{\gamma}} ,
[\bar{p} ]_q \bigg)  \quad \mbox{\rm and} \quad
 \delta_2 =\bigg( \frac{1}{q\sqrt{\gamma}} ,[\bar{p} ]_q
,0,[0]_q \bigg) ,$$
where $ \bar{p} \in \Z$ is such that $p_0 \bar{p} =-1
\mod{q}$.
Set  $ V_j =\pi_{\delta_j}^* =\pi_{-\delta_j}$, $j=1,2$. For any
$\phi \in {\mathcal S} (\R)$, we
 consider $\phi_1,\phi_2 \in {\mathcal S}(M)$ defined by
\begin{equation*}
\begin{split}
\phi_1 \big( s,[n]_q \big) & =\phi(s) , \\
 \phi_2 \big( s,[n]_q \big)&  =\sqrt{q}\,  \delta_{[0]_q , [n]_q } \phi(s)
 ,\quad s\in \R  ,[n]_q \in \Z_q  ,
\end{split}
\end{equation*}
where $\delta_{a,b}$, $a,b\in \Z_q$ denotes Kronecker's symbol.
Denote  also $\delta_a (b)=\delta_{a,b}$, $ a,b\in \Z_q$.
Notice that if ${\mathcal F} \phi =\phi$ on $\R$, then ${\mathcal F}
 (\phi_1 +\phi_2)=
\phi_1 +\phi_2$ on $M$ because ${\mathcal F} \big( 1+\sqrt{q}\, \delta_{[0]_q} \big)
=1+\sqrt{q}\, \delta_{[0]_q}$ (again, it is essential that the measure on $
\Z_q$ is the Haar-Plancherel one). If we set
$$b_{m_1,m_2} =\int\limits_{\R} \overline{\phi(s)} \,
\phi \bigg( s+\frac{m_2}{q\sqrt{\gamma}}\bigg)  {\mathbf e} \bigg(
\frac{sm_1}{q\sqrt{\gamma}}\bigg)  ds ,\quad m_1,m_2 \in \Z , $$
then
\begin{equation}
\begin{split}
 \langle \phi_1,\phi_1 \rangle_{D^\bot} (m_1 \delta_1 +m_2 \delta_2) &  =
\langle \phi_1,\phi_1 \rangle_{D^\bot} \bigg( \frac{m_2}{q\sqrt{\gamma}} ,
[m_2 \bar{p} ]_q  ,\frac{m_1}{q\sqrt{\gamma}}
  ,[m_1 \bar{p} ]_q \bigg) \\
 &  =\int\limits_{\R \times \Z_q} \overline{\phi_1 (s,[n]_q)}
\,
\phi_1 \bigg( s+\frac{m_2}{q\sqrt{\gamma}}  ,[n+m_2 \bar{p} ]_q \bigg)
{\mathbf e} \bigg(\frac{sm_1}{q\sqrt{\gamma}} +\frac{nm_1 \bar{p}}{q}\bigg)
ds  d[n]_q \\  &
=\frac{1}{\sqrt{q}} \sum\limits_{n\in \Z_q}
{\mathbf e} \bigg( \frac{nm_1 \bar{p}}{q}\bigg)
 b_{m_1,m_2}
=\begin{cases} 0, & \mbox{\rm if $q\nmid m_1 $} \\
\sqrt{q}\, b_{m_1,m_2}, & \mbox{\rm if $q\mid m_1 $} \end{cases}  ,
  \\
\langle \phi_1,\phi_2 \rangle_{D^\bot} (m_1 \delta_1 +m_2 \delta_2)  &
 =\int\limits_{\R \times \Z_q} \overline{\phi_1 (s,[n]_q)}\,
\phi_2 \bigg( s+\frac{m_2}{q\sqrt{\gamma}} \, ,[n+m_2 \bar{p} ]_q \bigg) \\ &  \qquad \times
{\mathbf e} \bigg( \frac{sm_1}{q\sqrt{\gamma}} +\frac{nm_1 \bar{p}}{q}\bigg)
ds  d[n]_q
 ={\mathbf e} \bigg( -\frac{m_1 m_2 \bar{p}^2}{q}\bigg) b_{m_1,m_2} \, ,\\
\langle \phi_2,\phi_1 \rangle_{D^\bot} (m_1 \delta_1 +m_2 \delta_2)  &
=b_{m_1,m_2} ,\\
\langle \phi_2,\phi_2 \rangle_{D^\bot} (m_1 \delta_1 +m_2 \delta_2) &
=\begin{cases} 0, & \mbox{\rm if $q \nmid  m_2$} \\
\sqrt{q}\, b_{m_1,m_2}, & \mbox{\rm if $q\mid m_2$} \end{cases} ,
\end{split}
\end{equation}
and consequently for all $m_1,m_2 \in \Z$,
$$\langle \phi_1 +\phi_2,\phi_1 +\phi_2 \rangle_{D^\bot} (m_1 \delta_1 +
m_2 \delta_2) =b_{m_1,m_2}\, c_{m_1,m_2}\, ,$$
where
$$c_{m_1,m_2} =\begin{cases} 1+{\mathbf e} \Big( -\frac{m_1 m_2 \bar{p}^2}{q}\Big), &
\mbox{\rm if $ q\nmid m_1$ and $q \nmid m_2$} \\
2+\sqrt{q}, & \mbox{\rm if $q \nmid m_1 ,\ q\mid m_2$ or $
q\mid m_1 ,\ q\nmid m_2$} \\
2+2\sqrt{q}, & \mbox{\rm if $q\mid m_1$ and $q\mid m_2$} \end{cases}  .$$

If $\phi(s)=e^{-\pi s^2}$, set $\tilde{\phi} =\phi_1 +\phi_2$.
Computations similar to those from Section 1 (with $q\sqrt{\gamma}$ instead of
$\sqrt{\alpha})$ yield for all $m_1,m_2 \in \Z$,
$$b_{m_1,m_2} =\frac{1}{\sqrt{2}} \, e^{-\frac{\pi (m_1^2+m_2^2)}{
2q^2\gamma} -\frac{\pi im_1 m_2}{q^2 \gamma} } ,$$
and
\begin{equation*}
\begin{split}
 a\sqrt{2} & =\sqrt{2} \langle \tilde{\phi} ,\tilde{\phi} \rangle_{D^\bot}
=\sum\limits_{m_1,m_2} b_{m_1,m_2} c_{m_1,m_2} \pi^*_{m_1\delta_1 +m_2\delta_2}
\\  &
 =\sqrt{2} \sum\limits_{m_1,m_2} b_{-m_1,-m_2}\ c_{-m_1,-m_2}
\beta (-m_2 \delta_2 ,-m_1\delta_1)  \pi_{m_1 \delta_1 +m_2 \delta_2}
 \\ &
=\sqrt{2} \sum\limits_{m_1,m_2} b_{-m_1,-m_2} c_{-m_1,-m_2}  V_2^{m_2}
 V_1^{m_1}
 =\sqrt{2}\sum\limits_{m_1,m_2} b_{m_1,m_2} c_{m_1,m_2}  V_2^{m_2} V_1^{m_1}
\\   &
 =\sum\limits_{m_1,m_2} c_{m_1,m_2}
e^{-\frac{\pi(m_1^2+m_2^2)}{2q^2 \gamma} -\frac{\pi im_1 m_2}{q^2 \gamma}}
 V_2^{m_2} V_1^{m_1} \\ &
 =\sum\limits_{m_2} e^{ -\frac{\pi m_2^2}{2q^2 \gamma}}
 \left( \sum\limits_{m_1} c_{m_1,m_2}
e^{\frac{2\pi im_1 m_2 \bar{p}^2}{q} +
  \frac{\pi im_1 m_2}{q^2 \gamma} -\frac{\pi m_1^2}{2q^2 \gamma}} V_1^{m_1} \right) V_2^{m_2} .
\end{split}
\end{equation*}
We set
$\alpha_m =e^{-\frac{\pi m^2}{2q^2 \gamma}}$, $m\in \Z$ and
$$D_n =\sum\limits_m c_{m,n} \alpha_m
{\mathbf e} \bigg( \frac{mn}{q} \bigg(
\frac{1}{2q\gamma} +\bar{p}^2 \bigg) \bigg) V_1^m ,\quad n\in \Z .$$
For all $f\in L^2 (\R \times \Z_q )$,
$$\big( V_1^m f\big) (x,[k]_q) ={\mathbf e} \bigg( -\frac{mx}{q\sqrt{\gamma}}
-\frac{mk\bar{p}}{q} \bigg) f(x,[k]_q)
,$$
hence
$D_n$ is the multiplication operator $M_{\Phi_n}$ on
 $L^2 (\R \times \Z_q)$, with $
 \Phi_n \in L^\infty (\R \times \Z_q)$
given by
$$\Phi_n (x,[k]_q) =\sum\limits_m c_{m,n}  \alpha_m
{\mathbf e} \bigg(\frac{mn}{q}
\bigg( \frac{1}{2q\gamma} +\bar{p}^2 \bigg) -
\frac{mx}{q\sqrt{\gamma}} -\frac{mk\bar{p}}{q} \bigg)$$
and $ a\sqrt{2} =\sum\limits_n \alpha_n D_n V_2^n$. We have
\begin{equation}
\begin{split}
 \Phi_0 (x,[k]_q) &
=\sum\limits_m c_{m,0}  \alpha_m
{\mathbf e} \bigg( -\frac{mx}{q\sqrt{\gamma}} -\frac{mk\bar{p}}{q} \bigg)  \\ &
 =(2+2\sqrt{q} ) \sum\limits_l \alpha_{ql}
{\mathbf e} \bigg(-\frac{lx}{\sqrt{\gamma}}\bigg) +
(2+\sqrt{q})\sum\limits_{q\,\nmid\,  m} \alpha_m
{\mathbf e} \bigg( -\frac{mx}{q\sqrt{\gamma}} -\frac{mk\bar{p}}{q}\bigg)  \\
&
 =\sqrt{q} \, \vartheta \bigg( -\frac{x}{\sqrt{\gamma}} \, ,\frac{i}{2\gamma} \bigg)
+(2+\sqrt{q} ) \vartheta \bigg( -\frac{x}{q\sqrt{\gamma}} -\frac{k\bar{p}}{q},
\frac{i}{2q^2 \gamma} \bigg) \\ &
 \geq (2+\sqrt{q} ) c\bigg( \frac{1}{2q^2 \gamma} \bigg) +\sqrt{q} \,
c\bigg( \frac{1}{2\gamma} \bigg) >0 ,
\end{split}
\end{equation}
where $c(t)=\min\limits_{x\in \R}  \, \vartheta (x,it)=
\vartheta \big(  \frac{1}{2}  ,it\big)$, $t>0$
and
\begin{equation}
\begin{split}
 \big| \Phi_n (x,[n]_q )\big|  & \leq \sum\limits_m \vert c_{m,n}
 \vert   \alpha_m  =\sum\limits_l c_{ql,n}  \alpha_{ql}
+\sum\limits_{q\, \nmid\, m} \vert c_{m,n} \vert
\alpha_m \\ &
 \leq (2+2\sqrt{q}) \vartheta \bigg( 0 ,\frac{i}{2\gamma} \bigg)
+(2+\sqrt{q} ) \bigg(  \vartheta \bigg( 0  ,
\frac{i}{2q^2\gamma} \bigg) -\vartheta \bigg( 0 ,
\frac{i}{2\gamma} \bigg) \bigg) \\  &
 =(2+\sqrt{q} ) C\bigg( \frac{1}{2q^2 \gamma} \bigg) +\sqrt{q}\,
C\bigg( \frac{1}{2\gamma} \bigg) ,
\end{split}
\end{equation}
where $C(t)=\max\limits_{x\in \R}  \,
\vartheta (x,it)=\vartheta (0,it)$, $t>0$. We combine (2.2) and (2.3) to
obtain for all $n\in \Z$,
\begin{equation}
 \bigg\| \frac{\Phi_n}{\Phi_0} \bigg\|_\infty  =\sup
\left\{
 \bigg|
\frac{\Phi_n (x,[k]_q)}{\Phi_0 (x,[k]_q)} \bigg|
\,  ;\, (x,[k]_q )\in \R \times \Z_q \right\}
\leq
\frac{(2+\sqrt{q}) C\big( \frac{1}{2q^2 \gamma} \big) +\sqrt{q}\,
C\big( \frac{1}{2\gamma} \big)}{(2+\sqrt{q})  c
\big( \frac{1}{2q^2 \gamma} \big) +\sqrt{q}\,
c\big( \frac{1}{2\gamma} \big)} .
\end{equation}

The function $h(t)$ from (1.28) is monotonically decreasing on $(0,\infty)$
and $C(t)<h(t) c(t)$ for all $t>0$, therefore (2.4) yields
for all $n\in \Z$,
$$\bigg\| \frac{\Phi_n}{\Phi_0} \bigg\|_\infty  \leq \max \left\{
h\bigg( \frac{1}{2q^2 \gamma} \bigg)\, ,\, h\bigg( \frac{1}{2\gamma} \bigg)
\right\} =h\bigg( \frac{1}{2q^2 \gamma} \bigg) ,$$
and
$$S=\sum\limits_{n\neq 0} \alpha_n \bigg\| \frac{\Phi_n}{\Phi_0} \bigg\|_\infty
\leq 2 h\bigg( \frac{1}{2q^2 \gamma} \bigg) \sum\limits_{n\geq 1}
\alpha_n .$$

But
$$\sum\limits_{n\geq 1} \alpha_n =\sum\limits_{n\geq 1} e^{-\frac{\pi n^2}{
2q^2 \gamma}} <\sum\limits_{n\geq 0} e^{-\frac{\pi (3n+1)}{2q^2 \gamma}}
=\frac{e^{-\frac{\pi}{2q^2\gamma}}}{1-e^{-\frac{3\pi}{2q^2 \gamma}}}  ,$$
therefore
$$S\leq 2 h\bigg( \frac{1}{2q^2 \gamma} \bigg) \frac{
e^{-\frac{\pi}{2q^2\gamma}}}{1-e^{-\frac{3\pi}{2q^2 \gamma}}} =g(r),$$
where $r=e^{-\frac{\pi}{2q^2 \gamma}}$ and $g$ is as in the proof of
Proposition 1.7. We have already shown there that $g(r)<1$ if
$\frac{1}{2q^2 \gamma} >0.531$, hence in this case
$I+\sum\limits_{n\neq 0} \alpha_n D_0^{-1} D_n V_2^n$ is bounded and invertible
on $L^2 (\R\times \Z_q )$, thus
$\sqrt{2} \langle \tilde{\phi} ,\tilde{\phi} \rangle_{D^\bot}
=\sum\limits_n \alpha_n D_n V_2^n$ is invertible.
Since the analogue of Proposition 1.1 with $\R$ replaced by $\R \times \Z_q$ holds, we get the following

\begin{prop}
Let $\alpha \in (0,1)$. If $q\geq 2$ and $p$ are integers such that
$q\geq 2$,
$0<\alpha - \frac{p}{q} <\frac{0.942}{q^2}$ and
there exists $p_0 \in \Z$ such that $p=p_0^2 \mod{q}$,
 then $a$ is invertible and
$A_\alpha^\sigma$ contains a projection of trace $q\alpha -p$.
\end{prop}

One proves in a similar way that if $0<\gamma^\prime = \frac{p}{q} -\alpha  <
\frac{0.942}{q^2}$ and there exists $p_1 \in \Z$ such that $p=-p_1^2
\mod{q}$, then $a$ is invertible, hence $A_\alpha^\sigma$ contains a projection
of trace $p-q\alpha$. The only difference is that one starts with
$\varepsilon_1 =\big( 0,[0]_q ,\sqrt{\gamma^\prime} ,[p_1 ]_q \big)$
and $\varepsilon_2 =\big( \sqrt{\gamma^\prime} ,[p_1 ]_q ,0,[0]_q \big)$.

\setcounter{equation}{0}

\section{More on the case $\alpha =q^{-1}$, $q\in \Z$, $q\geq 2$}

In this section we will focus mainly on the case $\alpha =
q^{-1}$, $q\in \Z$, $q\geq 2$, obtaining lower bounds
for the norm of $H_\alpha =\pi_{\varepsilon_1} +\pi_{\varepsilon_1}^* +
\pi_{\varepsilon_2} +\pi_{\varepsilon_2}^*$ and closed formulae for the projection
constructed in Corollary 1.5 (notation is as in Section 1).

In the beginning we will assume that $\alpha \in (0,1)$ is such that
the operator
$\langle f,f\rangle_{D^\bot}$ from (1.18) is invertible (for example
$0<\alpha <0.942$), hence
$$e=\left< f\langle f,f\rangle_{D^\bot}^{-\frac{1}{2}} ,
 f\langle f,f\rangle_{D^\bot}^{-\frac{1}{2}} \right> _D$$
is a projection of
trace $\alpha$ in $A_\alpha^\sigma$. Then, according to \cite{R1},
$\Theta :eA_\alpha e \simeq eC^* (D,\beta )e \rightarrow C^* (D^\bot ,
\overline{\beta} ) \simeq A_{\frac{1}{\alpha}}$ defined by
$$\Theta (exe)=\langle f,f\rangle_{D^\bot}^{-\frac{1}{2}}
\langle f,xf\rangle_{D^\bot} \langle f,f\rangle_{D^\bot}^{-\frac{1}{2}}  ,
\quad x\in C^* (D,\beta)$$
is a *-isomorphism. Direct computations based on (1.4) and (1.16) yield
\begin{equation*}
\begin{split}
 \langle f,\pi_{\pm \varepsilon_1} f\rangle_{D^\bot} & (m_1 \delta_1 +m_2
\delta_2 ) = \langle f,\pi_{\pm \varepsilon_1} f\rangle_{D^\bot} \bigg(
\frac{m_2}{\sqrt{\alpha}} ,\frac{m_1}{\sqrt{\alpha}} \bigg)
=\int\limits_{\R} \overline{f(s)} \, f\bigg( s+\frac{m_2}{\sqrt{\alpha}}
\pm \sqrt{\alpha} \bigg)   {\mathbf e} \bigg( \frac{sm_1}{\sqrt{\alpha}}\bigg)
 ds \\   &
=\int\limits_{\R} e^{-\pi s^2 -\pi \big( s+\frac{m_2}{\sqrt{\alpha}}\pm
\sqrt{\alpha} \big)^2} {\mathbf e} \bigg(\frac{sm_1}{\sqrt{\alpha}}\bigg)  ds
 =\frac{1}{\sqrt{2}} \, e^{-\frac{\pi \alpha}{2} -\frac{\pi (m_1^2 +m_2^2)}{2\alpha}
-\frac{\pi im_1 m_2}{\alpha} \mp \pi m_2 +\pi im_1}
\end{split}
\end{equation*}
and
\begin{equation*}
\begin{split}
 \langle f,\pi_{\pm \varepsilon_2} f\rangle_{D^\bot} (m_1 \delta_1 +m_2
\delta_2 ) &
=\int\limits_{\R} \overline{f(s)} \, f\bigg( s+\frac{m_2}{\sqrt{\alpha}}
\bigg)  {\mathbf e} \bigg( \pm \sqrt{\alpha} \bigg(
s+\frac{m_2}{\sqrt{\alpha}} \bigg)+\frac{sm_1}{\sqrt{\alpha}} \bigg)
 ds \\ & =\frac{1}{\sqrt{2}} \,
 e^{-\frac{\pi \alpha}{2} -\frac{\pi(m_1^2+m_2^2)}{2\alpha}
-\frac{\pi im_1 m_2}{\alpha} \mp \pi m_1 + \pi im_2} ,\end{split}
\end{equation*}
hence
\begin{equation}
\begin{split}
 \langle f,H_\alpha f\rangle_{D^\bot}  =\frac{1}{\sqrt{2}} &
\sum\limits_{m_1,m_2} e^{-\frac{\pi \alpha}{2} -
\frac{\pi(m_1^2+m_2^2)}{2\alpha}+\frac{\pi im_1 m_2}{\alpha}}  \\ &
 \times \left( (e^{\pi m_2} +e^{-\pi m_2} )e^{\pi im_1} +
(e^{\pi m_1} +e^{-\pi m_1} )e^{\pi im_2} \right) V_1^{m_1} V_2^{m_2} .
\end{split}
\end{equation}

If $\alpha =q^{-1},$ $q\in \Z$, $q\geq 2$, then
$A_{\frac{1}{\alpha}} \simeq C^* (D^\bot ,\overline{\beta} )$ is isomorphic
to $C(\T^2)$ and $\Theta (eH_\alpha e)$ is identified with the function
$$\Theta (eH_\alpha e)(z_1 ,z_2 )=\frac{F(z_1 ,z_2)}{G(z_1 ,z_2)} ,
\quad z_1 ,z_2 \in \T ,$$
where
\begin{equation*}
\begin{split}
& F(z_1 ,z_2)  =\sum\limits_{m_1 ,m_2} e^{-\frac{\pi \alpha}{2} -
\frac{\pi(m_1^2+m_2^2)}{2\alpha}+\frac{\pi im_1 m_2}{\alpha}}
\left( (e^{\pi m_2} +e^{-\pi m_2} )e^{\pi im_1} +
(e^{\pi m_1} +e^{-\pi m_1} )e^{\pi im_2} \right) z_1^{m_1} z_2^{m_2}  , \\ &
 G(z_1 ,z_2 ) =\sum\limits_{m_1,m_2} e^{-\frac{\pi(m_1^2+m_2^2)}{2\alpha} +
\frac{\pi i m_1 m_2}{\alpha}} z_1^{m_1} z_2^{m_2}  .
\end{split}
\end{equation*}
Since
$$\| H_\alpha \| \geq \| eH_\alpha e\| =\| \Theta (eH_\alpha e)\|_\infty
\geq \frac{F(1,1)}{G(1,1)}  ,$$
we further get for $q$ even
\begin{equation}
 \| H_{\frac{1}{q}} \| \geq \phi_0 \bigg( \frac{1}{q} \bigg)  =
4\, e^{-\frac{\pi}{2q}} \, \frac{\vartheta \big( \frac{i}{2},\frac{iq}{2} \big) \,
\vartheta \big( \frac{1}{2} ,\frac{iq}{2} \big)}{\vartheta \big( 0,
\frac{iq}{2} \big)^2}
 =4\, \frac{\vartheta \big( \frac{1}{q} ,\frac{2i}{q} \big) \, \vartheta \big(
\frac{1}{2} ,\frac{iq}{2} \big)}{ \vartheta \big( 0,\frac{2i}{q} \big)\,
\vartheta \big( 0 , \frac{iq}{2} \big)},
\end{equation}
and for $q$ odd
\begin{equation}
\| H_{\frac{1}{q}} \| \geq \phi_1  \bigg( \frac{1}{q} \bigg)  =
4\, e^{-\frac{\pi}{2q}} \, \frac{\vartheta \big( \frac{i}{2},\frac{iq}{2} \big)
\vartheta \big( \frac{1}{2} ,\frac{iq}{2} \big) -2\,
\vartheta^{\operatorname{odd}} \big( \frac{i}{2},\frac{iq}{2} \big)
\vartheta^{\operatorname{odd}} \big( \frac{1}{2} ,\frac{iq}{2} \big)}{
\vartheta \big( 0,
\frac{iq}{2} \big)^2 -2\, \vartheta^{\operatorname{odd}} \big( 0,
\frac{iq}{2} \big)^2} ,
\end{equation}
where we set
$$\vartheta^{\operatorname{odd}} (z,\tau)=\sum\limits_{m \operatorname{odd}}
e^{\pi im^2 \tau +2\pi imz} =\vartheta (z,\tau) -\vartheta (2z,4\tau ), \quad
z\in \C, \ \tau \in \HH .$$

Taking as usual
\begin{equation}
\vartheta_{a,b} (z,\tau) =\sum\limits_m
e^{\pi i (m+a)^2\tau +2\pi i(m+a)(z+b)} ,\quad  z\in \C,\ \tau \in
\HH ,\
a,b\in \R ,
\end{equation}
we make use of $\vartheta \big( \frac{i}{2} ,it \big) =
e^{\frac{\pi}{4t}} \vartheta_{\frac{1}{t},0} (0,it)$ and
$ \vartheta \big( \frac{1}{2}, it\big) =\vartheta_{0,\frac{1}{2}} (0,it )$ to get
$$\phi_0 (t)=\frac{4\vartheta_{0,\frac{1}{2}} \big( 0,\frac{i}{2t} \big)
\vartheta_{2t,0} \big( 0,\frac{i}{2t} \big)}{\vartheta \big( 0,\frac{i}{2t} \big)^2}
.$$

The graphs of $\phi_0$ and $\phi_1$ are drawn in Figure 3.1 and 3.2 using Mathematica.

\bigskip

\begin{figure}
\includegraphics[scale=0.55,bb= 100 10 250 250]{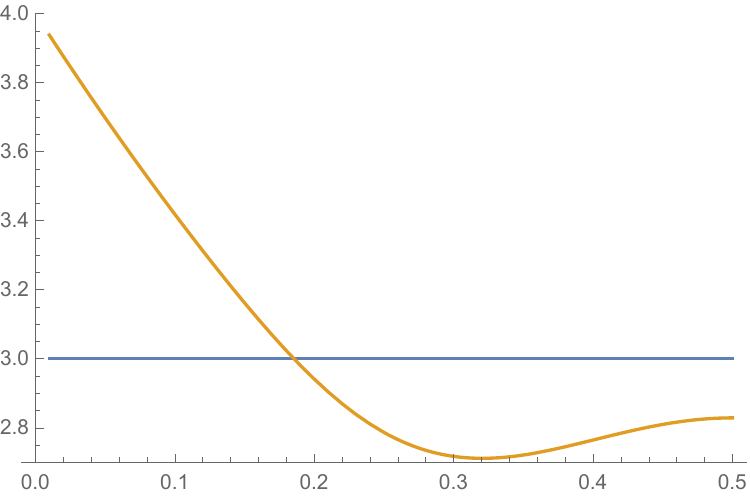}
\caption{The graph of $\phi_0$ on $\big( 0,\frac{1}{2}\big]$}\label{Figure 3.1}
\end{figure}

\begin{figure}
\includegraphics[scale=0.55,bb=100 10 250 250]{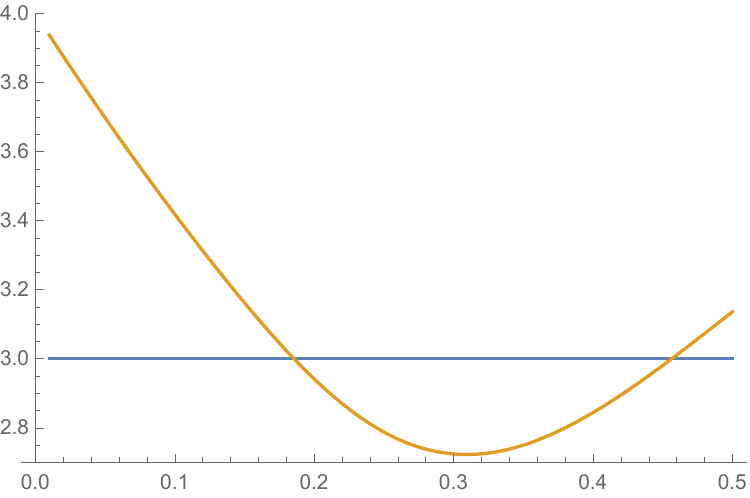}
\caption{The graph of $\phi_1$ on $\big( 0,\frac{1}{2}\big]$}\label{Figure 3.2}
\end{figure}

\bigskip

They should be compared with Hofstadter's butterfly (\cite{H}).

Estimates (3.2) and (3.3) are quite accurate, as suggested by the following
table. The norm of $H_{\frac{1}{q}}$ has been computed numerically using
\cite[Cor.3.2]{AMS}.

\bigskip

\hspace{1.5cm}
\begin{tabular}{|c|c|c|c|}
\hline \hline
\multicolumn{1}{|p{1.5cm}|}{\hspace{0.3cm} $\alpha$ }  &
\multicolumn{1}{|p{2.5cm}|}{\hspace{0.7cm} $\| H_\alpha \|$} &
\multicolumn{1}{|p{2.5cm}|}{\hspace{0.8cm}$\phi_0 (\alpha)$} &
\multicolumn{1}{|p{2.5cm}|}{\hspace{0.8cm}$\phi_1 (\alpha)$} \\
\hline \hline
\multicolumn{1}{|p{1cm}|}{\hspace{0.2cm} 1\slash 2} &
\multicolumn{1}{|p{2.5cm}|}{\hspace{0.7cm}2.82842} &
\multicolumn{1}{|p{2.5cm}|}{\hspace{0.7cm}2.82842} &
\multicolumn{1}{|p{2.5cm}|}{} \\
\hline
\multicolumn{1}{|p{1cm}|}{\hspace{0.2cm} 1\slash 3} &
\multicolumn{1}{|p{2.5cm}|}{\hspace{0.7cm}2.73205} &
\multicolumn{1}{|p{2.5cm}|}{\hspace{0.7cm}} &
\multicolumn{1}{|p{2.5cm}|}{\hspace{0.7cm}2.73205} \\
\hline
\multicolumn{1}{|p{1cm}|}{\hspace{0.2cm} 1\slash 4} &
\multicolumn{1}{|p{2.5cm}|}{\hspace{0.7cm}2.82842} &
\multicolumn{1}{|p{2.5cm}|}{\hspace{0.7cm}2.78648} &
\multicolumn{1}{|p{2.5cm}|}{\hspace{0.7cm}} \\
\hline
\multicolumn{1}{|p{1cm}|}{\hspace{0.2cm} 1\slash 5} &
\multicolumn{1}{|p{2.5cm}|}{\hspace{0.7cm}2.96645} &
\multicolumn{1}{|p{2.5cm}|}{\hspace{0.7cm}} &
\multicolumn{1}{|p{2.5cm}|}{\hspace{0.7cm}2.94109} \\
\hline
\multicolumn{1}{|p{1cm}|}{\hspace{0.2cm} 1\slash 6} &
\multicolumn{1}{|p{2.5cm}|}{\hspace{0.7cm}3.09557} &
\multicolumn{1}{|p{2.5cm}|}{\hspace{0.7cm}3.08292} &
\multicolumn{1}{|p{2.5cm}|}{\hspace{0.7cm}} \\
\hline
\multicolumn{1}{|p{1cm}|}{\hspace{0.2cm} 1\slash 7} &
\multicolumn{1}{|p{2.5cm}|}{\hspace{0.7cm}3.20330} &
\multicolumn{1}{|p{2.5cm}|}{\hspace{0.7cm}} &
\multicolumn{1}{|p{2.5cm}|}{\hspace{0.7cm}3.19690} \\
\hline
\multicolumn{1}{|p{1cm}|}{\hspace{0.2cm} 1\slash 8} &
\multicolumn{1}{|p{2.5cm}|}{\hspace{0.7cm}3.29066} &
\multicolumn{1}{|p{2.5cm}|}{\hspace{0.7cm}3.28709} &
\multicolumn{1}{|p{2.5cm}|}{\hspace{0.7cm}} \\
\hline
\multicolumn{1}{|p{1cm}|}{\hspace{0.2cm} 1\slash 9} &
\multicolumn{1}{|p{2.5cm}|}{\hspace{0.7cm}3.36165} &
\multicolumn{1}{|p{2.5cm}|}{\hspace{0.7cm}} &
\multicolumn{1}{|p{2.5cm}|}{\hspace{0.7cm}3.35943} \\
\hline
\multicolumn{1}{|p{1cm}|}{\hspace{0.2cm} 1\slash 10} &
\multicolumn{1}{|p{2.5cm}|}{\hspace{0.7cm}3.42005} &
\multicolumn{1}{|p{2.5cm}|}{\hspace{0.7cm}3.41855} &
\multicolumn{1}{|p{2.5cm}|}{\hspace{0.7cm}} \\
\hline
\multicolumn{1}{|p{1cm}|}{\hspace{0.2cm} 1\slash 11} &
\multicolumn{1}{|p{2.5cm}|}{\hspace{0.7cm}3.46880} &
\multicolumn{1}{|p{2.5cm}|}{\hspace{0.7cm}} &
\multicolumn{1}{|p{2.5cm}|}{\hspace{0.7cm}3.46771} \\
\hline
\multicolumn{1}{|p{1cm}|}{\hspace{0.2cm} 1\slash 12} &
\multicolumn{1}{|p{2.5cm}|}{\hspace{0.7cm}3.51004} &
\multicolumn{1}{|p{2.5cm}|}{\hspace{0.7cm}3.50922} &
\multicolumn{1}{|p{2.5cm}|}{\hspace{0.7cm}} \\
\hline
\multicolumn{1}{|p{1cm}|}{\hspace{0.2cm} 1\slash 13} &
\multicolumn{1}{|p{2.5cm}|}{\hspace{0.7cm}3.54537} &
\multicolumn{1}{|p{2.5cm}|}{\hspace{0.7cm}} &
\multicolumn{1}{|p{2.5cm}|}{\hspace{0.7cm}3.54473} \\
\hline
\multicolumn{1}{|p{1cm}|}{\hspace{0.2cm} 1\slash 50} &
\multicolumn{1}{|p{2.5cm}|}{\hspace{0.7cm}3.87630} &
\multicolumn{1}{|p{2.5cm}|}{\hspace{0.7cm}3.87628} &
\multicolumn{1}{|p{2.5cm}|}{\hspace{0.7cm}} \\
\hline
\multicolumn{1}{|p{1cm}|}{\hspace{0.2cm} 1\slash 51} &
\multicolumn{1}{|p{2.5cm}|}{\hspace{0.7cm}3.87869} &
\multicolumn{1}{|p{2.5cm}|}{\hspace{0.7cm}} &
\multicolumn{1}{|p{2.5cm}|}{\hspace{0.7cm}3.87867} \\
\hline
\multicolumn{1}{|p{1.5cm}|}{\hspace{0.2cm} 1\slash 100} &
\multicolumn{1}{|p{2.5cm}|}{\hspace{0.7cm}3.93766} &
\multicolumn{1}{|p{2.5cm}|}{\hspace{0.7cm}3.93765} &
\multicolumn{1}{|p{2.5cm}|}{\hspace{0.7cm}} \\
\hline
\multicolumn{1}{|p{1.5cm}|}{\hspace{0.2cm} 1\slash 101} &
\multicolumn{1}{|p{2.5cm}|}{\hspace{0.7cm}3.93827} &
\multicolumn{1}{|p{2.5cm}|}{\hspace{0.7cm}} &
\multicolumn{1}{|p{2.5cm}|}{\hspace{0.7cm}3.93827} \\
\hline \hline
\end{tabular}

\bigskip

We notice equality holds in (3.2) if $q=2$. To see this,
remark first that $H^2_{\frac{1}{2}} =4+H_0$, which yields
$\| H_{\frac{1}{2}} \| =2\sqrt{2}$. Using (3.4) and (1.15), we readily
see that
$$\vartheta_{0,\frac{1}{2}} (0,i)=\vartheta \bigg( \frac{1}{2}\, ,i\bigg) =
\vartheta_{\frac{1}{2},0} (0,i) .$$

On the other hand, Jacobi's identity (\cite[p.23]{M}) yields
$$\vartheta (0,i)^4 =\vartheta_{0,\frac{1}{2}} (0,i)^4 +
\vartheta_{\frac{1}{2},0} (0,i)^4 ,$$
hence $\vartheta \big( \frac{1}{2} ,i\big) =2^{-\frac{1}{4}}
\vartheta (0,i)$ and we get equality in (3.2) in this case.


Notice also that the numerical computations above suggest that
equality holds eventually for $\frac{1}{3}$ too, which would produce
an interesting relation between $\vartheta \big(  \frac{1}{2} ,
\frac{3i}{2} \big)$, $\vartheta \big( \frac{1}{3} ,\frac{2i}{3} \big)$,
$\vartheta \big( \frac{1}{6} ,\frac{i}{6} \big)$ and
$\vartheta (0,6i)$.

Next, we derive explicit formulae for the projection $e$ of trace $\alpha$
 constructed in Proposition 1.1 when $\alpha = q^{-1} $,
$q\in \Z$, $q\geq 2$. Again, notation is as in Section 1.

For any unitary operator $U$ and $\tau \in \HH$, $b\in \R$,
 $N\in {\mathbf N}^*$,
$a\in \Z$, we define the bounded operators
\begin{equation*}
\begin{split}
\vartheta (U,\tau) & =\sum\limits_m e^{\pi i\tau m^2} U^m ,\\
 \vartheta_{0,b} (U,\tau ) & =\sum\limits_m e^{\pi i\tau m^2 +2\pi imb} U^m , \\
 \vartheta_{\frac{a}{N},b}^{(q)} (U,\tau) & =\sum\limits_m
e^{\pi i\tau ( m+\frac{a}{N} )^2 +2\pi i( m+\frac{a}{N} ) b}
U^{Nm+a} .
\end{split}
\end{equation*}

The explicit expression of $e$ is given by the following

\begin{prop}
{\em (i)} If $q$ is even, then
$$e=
\frac{\displaystyle \sum\limits_{r,s=0}^{q-1}e^{-\frac{\pi irs}{q}} \vartheta^{(q)}_{\frac{s}{q},\frac{r}{2}} \big( U_2 ,\tfrac{iq}{2} \big)
\vartheta^{(q)}_{\frac{r}{q},\frac{s}{2}} \big( U_1 ,\tfrac{iq}{2} \big)}{q
\vartheta \big( U_2^q ,\frac{iq}{2} \big)
\vartheta \big( U_1^q ,\frac{iq}{2} \big)} .$$
{\em (ii)} If $q$ is odd, then
$$ e =
\frac{\displaystyle \sum\limits_{r,s=0}^{2q-1} e^{\frac{\pi irs}{q}+\pi i [ \frac{r}{q}][ \frac{s}{q}]}
 \vartheta^{(2q)}_{\frac{s}{2q} ,0} (U_2 ,2iq) \vartheta^{(2q)}_{\frac{r}{2q},0} (U_1 ,2iq)
}{q\displaystyle \sum\limits_{\varepsilon_1,\varepsilon_2 =0}^1 e^{
\pi i\varepsilon_1 \varepsilon_2} \vartheta^{(2)}_{\frac{\varepsilon_2}{2},0}
(U_2 ,2iq) \vartheta^{(2)}_{\frac{\varepsilon_1}{2},0} (U_1, 2iq)}  .
 $$
\end{prop}

{\bf Remark}.
The denominators are central elements in $A_{\frac{1}{q}}$.

Any $b\in C^* (D^\bot ,\overline{\beta}) =A_{\frac{1}{\alpha}}$ is represented as
\begin{equation}
b=\sum\limits_{z\in D^\bot} b(z)\pi_z^* =\sum\limits_{m_1,m_2}
b_{m_1,m_2} \pi^*_{m_1 \delta_1 +m_2\delta_2} =\sum\limits_{m_1,m_2}
b_{-m_1,-m_2} {\mathbf e} \bigg( \frac{m_1 m_2}{\alpha}\bigg) V_1^{m_1} V_2^{m_2} .
\end{equation}

For such $b$ and $n_1 ,n_2 \in \Z$, we set
\begin{equation}
\alpha_{n_1,n_2} =\langle fb,fb\rangle_D (n_1 \sqrt{\alpha} ,n_2 \sqrt{\alpha} ).
\end{equation}

\begin{lem}
Let $\alpha \in (0,1)$ and $b,\alpha_{n_1,n_2}$ be as in
 {\em (3.5)} and {\em (3.6)}. Then for all $n_1,n_2 \in \Z$,
\begin{equation*}
\begin{split}
\alpha_{n_1,n_2}  & =\frac{1}{\sqrt{2}} \,
e^{-\frac{\pi \alpha (n_1^2+n_2^2)}{2} +\pi in_1 n_2 \alpha} \\ &
\times \sum\limits_{m_1,\ldots,m_4} \Big(
 b_{m_1,m_2} \, \overline{b_{m_3,m_4}} \,
 e^{-\frac{\pi ( (m_1-m_3)^2 +(m_2-m_4)^2 )}{2\alpha} -\pi i
(n_1 +in_2)( m_3 -m_1 +i(m_4 -m_2)) +\frac{\pi i(m_1+m_3)(m_2 -m_4)}{\alpha}}
\Big) .
\end{split}
\end{equation*}
\end{lem}

\begin{proof}
Using
$$(fb)(s)=\sum\limits_{m_1,m_2} b_{m_1,m_2} e^{\frac{2\pi im_1 m_2}{\alpha}
-\pi ( s-\frac{m_2}{\sqrt{\alpha}} )^2 -\frac{2\pi ism_1}{\sqrt{\alpha}}}
,\quad s\in \R$$
and (1.3) we get
\begin{equation*}
\begin{split}
\alpha_{n_1,n_2} &  =
\int\limits_{\R} (fb)(s) \, \overline{(fb)(s+n_1 \sqrt{\alpha})} \,
e^{-2\pi isn_2 \sqrt{\alpha}} ds \\ &
 =\sum\limits_{m_1,\ldots,m_4} b_{m_1,m_2} \, \overline{b_{m_3,m_4}}\,
e^{\frac{2\pi i(m_1m_2 -m_3 m_4)}{\alpha} -\frac{\pi (m_2^2 +m_4^2)}{\alpha}
-\pi n_1^2 \alpha +2\pi n_1 m_4} \\ &
 \qquad
 \times \int\limits_{\R} e^{-2\pi s^2 -2\pi s \big(
-\frac{m_2+m_4}{\sqrt{\alpha}} +n_1 \sqrt{\alpha} +i(
\frac{m_1-m_3}{\sqrt{\alpha}} +n_2 \sqrt{\alpha} ) \big)} ds .
\end{split}
\end{equation*}
The statement follows now through a plain computation based on (1.16).
\end{proof}

For all $r,s,n_1,n_2 \in \Z$ we set
\begin{equation}
 a_{r,s}^{(n_1,n_2)} =e^{-\frac{\pi(r^2 +s^2)}{2\alpha} -\frac{\pi irs}{\alpha}
-\pi i(n_1+in_2)(r+is)},
\end{equation}
\begin{equation}
a^{(n_1,n_2)} =\sum\limits_{r,s} a_{r,s}^{(n_1,n_2)} V_1^r V_2^s  .
\end{equation}

Since
$$\sum\limits_{r,s} \vert a_{r,s}^{(n_1,n_2)} \vert  =
\vartheta \bigg( \frac{in_1}{2} ,\frac{i}{2\alpha} \bigg)
\vartheta \bigg( \frac{in_2}{2}  ,\frac{i}{2\alpha} \bigg)  <\infty  ,$$
it follows that $a^{(n_1,n_2)} \in C^* (D^\bot,\overline{\beta})$. Although we are not
going to use it, notice that since $\sigma_{D^\bot} (V_2)=V_1$ and
$\sigma_{D^\bot} (V_1)=V_2^{-1}$, we get
$$\sigma_{D^\bot} \big( a^{(n_1,n_2)} \big) =a^{(-n_2,n_1)} .$$

If $\alpha = q^{-1}$, $q\in \Z$, $q\geq 2$, Lemma 3.2 yields
for all $n_1,n_2 \in \Z$,
\begin{equation}
\alpha_{n_1,n_2} =\frac{1}{\sqrt{2}}\, e^{-\frac{\pi(n_1^2+n_2^2)}{2q}
+\frac{\pi in_1 n_2}{q}}
\sum\limits_{m_1,\ldots,m_4}  b_{m_1,m_2} \, \overline{b_{m_3,m_4}}\,
a_{m_3 -m_1 ,m_4 -m_2}^{(n_1,n_2)} .
\end{equation}

Moreover, for any $b$ as in (3.5),
\begin{equation*}
\begin{split}
  \tau_{D^\bot} \big( ba^{(n_1,n_2)} b^* \big) &
=\sum\limits_{m_1,\ldots,m_4} \sum\limits_{r,s}  b_{m_1,m_2}
 \overline{b_{m_3,m_4}}\, a_{r,s}^{(n_1,n_2)}
 \tau_{D^\bot} (V_1^{m_1+r-m_3} V_2^{m_2+s-m_4} ) \\ &
 =\sum\limits_{m_1,\ldots,m_4}   b_{m_1,m_2} \,
 \overline{b_{m_3,m_4}}\, a^{(n_1,n_2)}_{m_3-m_1,m_4-m_2}  ,
\end{split}
\end{equation*}
which together with (3.9) gives
\begin{equation*}
\begin{split}
\alpha_{n_1,n_2} &  =
\frac{1}{\sqrt{2}} \, e^{-\frac{\pi(n_1^2+n_2^2)}{2q} +\frac{\pi in_1n_2}{q}}
\tau_{D^\bot} \big( ba^{(n_1,n_2)} b^* \big) \\ &
 =\frac{1}{\sqrt{2}} \, e^{-\frac{\pi(n_1^2+n_2^2)}{2q} +\frac{\pi in_1n_2}{q}}
\int\limits_{[0,1]^2} a^{(n_1,n_2)} (e^{2\pi it_1},e^{2\pi it_2} )\cdot
\vert b( e^{2\pi it_1},e^{2\pi it_2} ) \vert^2  dt_1  dt_2 .
\end{split}
\end{equation*}

In the sequel we will take $b=a^{-\frac{1}{2}} =\langle f,f\rangle_{D^\bot}^{-
\frac{1}{2}} \in C(\T^2)$, hence
\begin{equation}
\begin{split}
\alpha_{n_1,n_2} &  =
\left< f\langle f,f\rangle_{D^\bot}^{-\frac{1}{2}} ,
f\langle f,f\rangle_{D^\bot}^{-\frac{1}{2}} \right>_D (n_1 \sqrt{\alpha} ,
n_2 \sqrt{\alpha} )  \\  &
 =\frac{1}{\sqrt{2}} \, e^{-\frac{\pi(n_1^2+n_2^2)}{2q} +\frac{\pi in_1n_2}{q}}
\int\limits_{[0,1]^2}
\frac{ a^{(n_1,n_2)} (e^{2\pi it_1},e^{2\pi it_2} )}{
a(e^{2\pi it_1},e^{2\pi it_2} )}  dt_1  dt_2
\end{split}
\end{equation}
and
\begin{equation}
e=\vert G\slash D \vert \sum\limits_{n_1,n_2} \alpha_{n_1,n_2} U_2^{n_2} U_1^{n_1}
=\frac{1}{q} \sum\limits_{n_1,n_2} \alpha_{n_1,n_2} U_2^{n_2} U_1^{n_1}  .
\end{equation}

For simplicity we shall denote throughout this section
$$\vartheta (z)=\vartheta \bigg( z,\frac{iq}{2} \bigg)  .$$

\begin{proof}[Proof of Proposition {\em 3.1 (i)}]
If $q$ is even, then (1.18) yields for all $t_1,t_2 \in \R$,
\begin{equation}
a(e^{2\pi it_1},e^{2\pi it_2} ) =\frac{1}{\sqrt{2}}\,
\vartheta (t_1)\vartheta (t_2)  .
\end{equation}
We employ (3.7) and (3.8) to obtain
\begin{equation}
a^{(n_1,n_2)} (e^{2\pi it_1},e^{2\pi it_2} ) =
\vartheta \bigg( t_1 -\frac{n_1+in_2}{2} \bigg)  \vartheta \bigg( t_2 -
\frac{n_2+in_1}{2} \bigg) .
\end{equation}
Writing $n_1 =ql_1+r$, $n_2 =ql_2+s$, $l_1,l_2,r,s\in \Z$, $0\leq r,s<q$
and employing (3.10), (3.12) and (3.13) we get
$$
 \alpha_{n_1,n_2}  =
e^{-\frac{\pi q(l_1^2+l_2^2)}{2} -\pi (l_1r+l_2s)-\frac{\pi(r^2+s^2)}{2q} +
\pi i(l_1s+l_2r)+\frac{\pi irs}{q}}
 \int\limits_0^1 \frac{\vartheta \big( t_1 -\frac{r+is+il_2q}{2} \big)}{
\vartheta (t_1)}  dt_1
  \int\limits_0^1 \frac{\vartheta \big(
t_2 -\frac{s+ir+il_1 q}{2} \big)}{\vartheta (t_2)}  dt_2  .
$$
Using also the quasi-periodicity relation
$$\vartheta \bigg( z+\frac{ilq}{2} \bigg) =e^{\frac{\pi ql^2}{2} -2\pi ilz}
\vartheta (z), \quad  z\in \C,\ l\in \Z ,$$
we get
$$\alpha_{ql_1+r,ql_2+s} =e^{-\frac{\pi(r^2+s^2)}{2q}+\frac{\pi irs}{q}}
\int\limits_0^1 \frac{e^{2\pi il_2 t_1}  \vartheta \big( t_1 -
\frac{r+is}{2} \big)}{\vartheta (t_1)}  dt_1 \
\int\limits_0^1 \frac{e^{2\pi il_1 t_2} \vartheta \big( t_2 -
\frac{s+ir}{2} \big)}{\vartheta (t_2)}  dt_2 .$$
Furthermore, we make use of
$$\sum\limits_l e^{2\pi ily} \int\limits_0^1
\frac{e^{2\pi ilx} \, \vartheta \big( x-\frac{r+is}{2} \big)}{\vartheta (x)}  dx
=\sum\limits_l \int\limits_0^1 \frac{e^{2\pi ilx}
\vartheta \big( x-y-\frac{r+is}{2} \big)}{\vartheta (x-y)}  dx
 =\frac{\vartheta \big( y+\frac{r+is}{2} \big)}{\vartheta (y)}
$$
and of the fact that $U_1^q$ and $U_2^q$ are central in $A_{\frac{1}{\alpha}}$ to get
$$ \sum\limits_{l_1,l_2} \alpha_{ql_1+r,ql_2+s} U_2^{ql_2} U_1^{ql_1}
 =e^{-\frac{\pi(r^2+s^2)}{2q}+\frac{\pi irs}{q}}
\frac{\vartheta \big( e^{\pi i(r+is)}  U_2^q ,\frac{iq}{2} \big)
\vartheta \big( e^{\pi i(s+ir)}  U_1^q ,\frac{iq}{2} \big)}{
\vartheta \big( U_2^q ,\frac{iq}{2} \big)  \vartheta \big( U_1^q ,
\frac{iq}{2} \big)} .
$$

The equality from Proposition 3.1 (i) follows from this and from
\begin{equation}
e^{-\frac{\pi s^2}{2q} +\frac{\pi irs}{q}}  U^s  \vartheta \bigg(
e^{\pi i(r+is)}  U^q ,\frac{iq}{2} \bigg) =
\vartheta^{(q)}_{\frac{s}{q},\frac{r}{2}} \bigg( U,\frac{iq}{2} \bigg) .\qedhere
\end{equation}
\end{proof}

Before completing the proof of Proposition 3.1 for $q$ odd, we introduce
some notation by setting for all $z\in \C$, $\tau \in \HH$,
$a,b\in \R$,
\begin{equation*}
\begin{split}
 \vartheta_{a,b}^{\operatorname{odd}} (z,\tau) & =\sum\limits_{m \operatorname{odd}}
e^{\pi i\tau (m+a)^2 +2\pi i(m+a)(z+b)} =\vartheta_{a,b} (z,\tau)-
\vartheta_{\frac{a}{2},2b} (2z,4\tau) =\vartheta_{\frac{a+1}{2},2b} (2z,4\tau) , \\
&
\vartheta_{a,b}^{\operatorname{even}} =\sum\limits_{m \operatorname{even}}
e^{\pi i\tau (m+a)^2 +2\pi i(m+a)(z+b)} =\vartheta_{\frac{a}{2},2b} (2z,4\tau) , \\
&
\vartheta^{\diamond} (z,\tau)=\vartheta^{\diamond}_{0,0} (z,\tau) ,
\quad \diamond \in \{ \operatorname{even},\operatorname{odd}\} \ .
\end{split}
\end{equation*}
For all $z\in \C$, $\tau \in \HH$, $a,b\in \R$,
 $l\in \Z$  we have
\begin{equation}
 \vartheta^{\operatorname{even}} \bigg( z+\frac{1}{2} ,\tau \bigg) =\vartheta^{\operatorname{even}} (z,\tau)  ,
\end{equation}
\begin{equation}
 \vartheta^{\operatorname{odd}} \bigg( z+\frac{1}{2},\tau \bigg) =-\vartheta^{\operatorname{odd}} (z,\tau)  ,
\end{equation}
\begin{equation}
\vartheta_{a,b}^{\operatorname{odd}} (z+l\tau ,\tau )=
e^{-\pi i\tau l^2 -2\pi il(z+b)} \cdot \begin{cases}
\vartheta_{a,b}^{\operatorname{odd}} (z,\tau ), & \mbox{\rm if $l$ is even} \\
\vartheta_{a,b}^{\operatorname{even}} (z,\tau ), & \mbox{\rm if $l$ is odd} \end{cases} \ .
\end{equation}

The following lemma is easy to prove.

\begin{lem}
Let $\phi :\R^2 \rightarrow \C$
 be a continuous function which is periodic
modulo $\Z^2$. Then
\begin{equation*}
\begin{split}
& \sum\limits_{l_1,l_2 \operatorname{even}} \
\int\limits_{[0,1]^2}  e^{2\pi i(l_1 t_1 +l_2 t_2 )} \phi (t_1,t_2) dt_1  dt_2
=\frac{1}{4} \bigg(\phi (0,0) +\phi \bigg( \frac{1}{2} ,0\bigg) +
\phi \bigg( 0,\frac{1}{2} \bigg) +\phi \bigg( \frac{1}{2} ,\frac{1}{2} \bigg) \bigg)
, \\  &
 \sum\limits_{l_1 \operatorname{even},l_2 \operatorname{odd}} \
\int\limits_{[0,1]^2}   e^{2\pi i(l_1 t_1 +l_2 t_2 )} \phi (t_1,t_2)  dt_1
 dt_2
=\frac{1}{4} \bigg( \phi (0,0) +\phi \bigg( \frac{1}{2}  ,0\bigg)
 -\phi \bigg( 0,\frac{1}{2} \bigg) -\phi \bigg( \frac{1}{2} ,\frac{1}{2} \bigg) \bigg)
,\\
&
 \sum\limits_{l_1 \operatorname{odd},l_2 \operatorname{even}}  \
\int\limits_{[0,1]^2}  e^{2\pi i(l_1 t_1 +l_2 t_2 )} \phi (t_1,t_2)
 dt_1  dt_2
=\frac{1}{4} \bigg( \phi (0,0) -\phi \bigg( \frac{1}{2} ,0\bigg)
 +\phi \bigg( 0,\frac{1}{2} \bigg) -\phi
\bigg( \frac{1}{2} ,\frac{1}{2} \bigg) \bigg)  , \\ &
 \sum\limits_{l_1,l_2  \operatorname{odd}}  \
\int\limits_{[0,1]^2}  e^{2\pi i(l_1 t_1 +l_2 t_2 )} \phi (t_1,t_2)
 dt_1  dt_2
=\frac{1}{4} \bigg( \phi (0,0) -\phi \bigg( \frac{1}{2} ,0\bigg)
 -\phi \bigg( 0,\frac{1}{2} \bigg) +\phi \bigg( \frac{1}{2},\frac{1}{2} \bigg) \bigg)
.\end{split}
\end{equation*}
\end{lem}

\begin{proof}
 We consider for instance the second equality. The others
are similar. We set
$l_1 =2m_1$, $l_2 =2m_2 +1$, $x_1 =2t_1$, $x_2 =2t_2$, and divide
$[0,2]^2$ into four equal squares to get:
\begin{equation*}
\begin{split}
S_{01} & =\sum\limits_{l_1 \operatorname{even},l_2 \operatorname{odd}}  \
\int\limits_{[0,1]^2}  e^{2\pi i(l_1 t_1 +l_2 t_2 )} \phi (t_1,t_2)
 dt_1   dt_2  \\  &  =
\frac{1}{4} \sum\limits_{m_1,m_2} \  \int\limits_{[0,2]^2}
e^{2\pi i(m_1 x_1 +m_2 x_2)}   e^{\pi ix_2} \, \phi \bigg( \frac{x_1}{2}   ,
\frac{x_2}{2} \bigg)  dx_1  dx_2 \\ &
 =\frac{1}{4}  \sum\limits_{m_1,m_2} \ \int\limits_{[0,1]^2}
e^{2\pi i(m_1 x_1 +m_2 x_2)}  \psi (x_1 ,x_2)  dx_1  dx_2  ,
\end{split}
\end{equation*}
where
$$\psi (x_1 ,x_2)
 =e^{\pi ix_2}  \phi \bigg( \frac{x_1}{2}  ,\frac{x_2}{2} \bigg)
+ e^{\pi ix_2}  \phi \bigg( \frac{x_1+1}{2}  ,\frac{x_2}{2} \bigg) -
e^{\pi ix_2}   \phi \bigg( \frac{x_1}{2}  ,\frac{x_2 +1}{2} \bigg)  -
e^{\pi ix_2}   \phi \bigg( \frac{x_1 +1}{2}  ,\frac{x_2+1}{2} \bigg)  .
$$

Since $\psi (x_1 +1,x_2)=\psi (x_1,x_2 )=\psi (x_1 ,x_2 +1)$, $x_1 ,x_2 \in \R$,
we get
\begin{equation*}
S_{01} =\frac{1}{4}  \psi (0,0)=\frac{1}{4}  \bigg(
\phi (0,0)+\phi \bigg( \frac{1}{2} ,0\bigg) -\phi \bigg( 0,\frac{1}{2} \bigg)
-\phi \bigg( \frac{1}{2} ,\frac{1}{2} \bigg) \bigg) .\qedhere
\end{equation*}
\end{proof}

For all $a,b\in \R$, $z\in \C$, we set
$$
 \vartheta_{a,b} (z)=\vartheta_{a,b} \bigg( z,\frac{iq}{2} \bigg)  ,\qquad
 \vartheta_{a,b}^\diamond (z) =\vartheta_{a,b}^\diamond \bigg( z,\frac{iq}{2} \bigg)
, \quad \diamond \in \{ \operatorname{even},\operatorname{odd} \}  , $$
and consider
 the following continuous functions on $\R^2 \slash \Z^2$:
\begin{equation*}
\begin{split}
& \phi_{r,s}^{0,0} (t_1,t_2)=\frac{\vartheta \big( t_1 +\frac{r+is}{2} \big)
\vartheta \big( t_2 +\frac{s+ir}{2} \big) -2
\vartheta^{\operatorname{odd}} \big( t_1 +\frac{r+is}{2} \big)
\vartheta^{\operatorname{odd}} \big( t_2 +\frac{s+ir}{2} \big)}{\vartheta (t_1) \vartheta (t_2) -2
\vartheta^{\operatorname{odd}} (t_1 ) \vartheta^{\operatorname{odd}} (t_2)}  ,   \\
& \phi_{r,s}^{0,1} (t_1,t_2) =\frac{\vartheta \big( t_1 +\frac{r+is}{2} \big)
\vartheta_{0,\frac{1}{2}} \big( t_2 +\frac{s+ir}{2} \big) +2
\vartheta^{\operatorname{even}} \big( t_1 +\frac{r+is}{2} \big)
\vartheta^{\operatorname{odd}} \big( t_2 +\frac{s+ir}{2} \big)}{\vartheta (t_1)  \vartheta (t_2) -2
\vartheta^{\operatorname{odd}} (t_1 ) \vartheta^{\operatorname{odd}} (t_2)}  , \\
& \phi_{r,s}^{1,0} (t_1,t_2) =\frac{\vartheta_{0,\frac{1}{2}}
\big( t_1 +\frac{r+is}{2} \big)
\vartheta \big( t_2 +\frac{s+ir}{2} \big) +2
\vartheta^{\operatorname{odd}} \big( t_1 +\frac{r+is}{2} \big)
\vartheta^{\operatorname{even}} \big( t_2 +\frac{s+ir}{2} \big)}{\vartheta (t_1) \vartheta (t_2) -2
\vartheta^{\operatorname{odd}} (t_1 ) \vartheta^{\operatorname{odd}} (t_2)}  , \\
& \phi_{r,s}^{1,1} (t_1,t_2) =\frac{\vartheta_{0,\frac{1}{2}}
 \big( t_1 +\frac{r+is}{2} \big)
\vartheta_{0,\frac{1}{2}} \big( t_2 +\frac{s+ir}{2} \big) -2
\vartheta^{\operatorname{even}} \big( t_1 +\frac{r+is}{2} \big)
\vartheta^{\operatorname{even}} \big( t_2 +\frac{s+ir}{2} \big)}{\vartheta (t_1) \vartheta (t_2) -2
\vartheta^{\operatorname{odd}} (t_1 ) \vartheta^{\operatorname{odd}} (t_2)}  .
\end{split}
\end{equation*}

\begin{proof}[Proof of Proposition {\em 3.1 (ii)}]
Our first aim is to compute $\sum\limits_{l_1,l_2} \alpha_{ql_1+r,
ql_2+s} U_2^{ql_2} U_1^{ql_1}$. This is split into the following four sums:
\begin{equation*}
\begin{split}
& S_{00} =\sum\limits_{l_1,l_2 \operatorname{even}} \alpha_{ql_1+r,ql_2+s}
e^{2\pi i(l_1 x_1+l_2x_2)}, \qquad
S_{01} =\sum\limits_{l_1 \operatorname{even},l_2 \operatorname{odd}} \alpha_{ql_1+r,ql_2+s}
e^{2\pi i(l_1 x_1+l_2x_2)} , \\  &
S_{10} =\sum\limits_{l_1 \operatorname{odd},l_2 \operatorname{even}} \alpha_{ql_1+r,ql_2+s}
e^{2\pi i(l_1 x_1+l_2x_2)} , \qquad
S_{11} =\sum\limits_{l_1,l_2 \operatorname{odd}}\alpha_{ql_1+r,ql_2+s}
e^{2\pi i(l_1 x_1+l_2x_2)}  .\end{split}
\end{equation*}
From (1.18) we get for all $t_1,t_2 \in R$,
\begin{equation}
a(e^{2\pi it_1} ,e^{2\pi it_2} ) =\frac{1}{\sqrt{2}}
\left( \vartheta (t_1)\, \vartheta (t_2) -2\, \vartheta^{\operatorname{odd}} (t_1)
\vartheta^{\operatorname{odd}} (t_2) \right) >0 ,
\end{equation}
whilst from (3.7) and (3.8)
\begin{equation}
\begin{split}
 a^{(n_1,n_2)} (e^{2\pi it_1} ,e^{2\pi it_2} ) & =\frac{1}{\sqrt{2}}
\left( \vartheta \bigg( t_1 -\frac{n_1+in_2}{2} \bigg)
\vartheta \bigg( t_2 -\frac{n_2+in_1}{2} \bigg) \right.  \\
&  \qquad -2
\vartheta^{\operatorname{odd}} \left. \bigg( t_1 -\frac{n_1+in_2}{2} \bigg)
\vartheta^{\operatorname{odd}} \bigg( t_2 -\frac{n_2+in_1}{2} \bigg) \right)  .
\end{split}
\end{equation}

To compute $S_{00}$, we notice that the quasi-periodicity properties
(3.15), (3.16), (3.17) and
\begin{equation}
\vartheta \bigg( z+\frac{ilq}{2} \bigg) =e^{\frac{\pi ql^2}{2} -2\pi ilz}
\vartheta (z)\, ,\quad z\in \C ,\ l\in \Z
\end{equation}
yield for all $t_1,t_2 \in \R$, $n_1,n_2 \in \Z$, $n_1=ql_1+r$,
$n_2 =ql_2+s$, $l_1,l_2,r,s\in \Z$, $0\leq r,s<q$,
\begin{equation}
\begin{split}
& a^{(n_1,n_2)} (e^{2\pi it_1} ,e^{2\pi it_2} )
=\frac{1}{\sqrt{2}} \, e^{\frac{\pi(l_1^2+l_2^2)q}{2} -2\pi i(l_2t_1+l_1t_2)
-\pi i(l_2r+l_1s)+\pi(l_2 s+l_1r)} \\ & \qquad
 \times \left( \vartheta \bigg( t_1 -\frac{r+is}{2} \bigg) \vartheta
\bigg( t_2 -\frac{s+ir}{2} \bigg)  -2 \vartheta^{\operatorname{odd}} \bigg(
t_1-\frac{r+is}{2} \bigg) \vartheta^{\operatorname{odd}} \bigg( t_2 -\frac{s+ir}{2} \bigg)
\right)  .
\end{split}
\end{equation}
We combine (3.10), (3.18), (3.21) and the definition of $\phi_{r,s}^{0,0}$ to get
$$\alpha_{n_1,n_2} =e^{-\frac{\pi(r^2+s^2)}{2q}+\frac{\pi irs}{q}}
\int\limits_{[0,1]^2} e^{2\pi i(l_2t_1+l_1t_2)}
\phi_{r,s}^{0,0} (-t_1,-t_2) dt_1  dt_2  .$$

Using Lemma 3.3 we gather
\begin{equation}
\begin{split}
 S_{00} & =e^{-\frac{\pi(r^2+s^2)}{2q}+\frac{\pi irs}{q}}
 \sum\limits_{l_1,l_2 \operatorname{even}} \ \int\limits_{[0,1]^2}
e^{2\pi il_2 (t_1+x_2)+2\pi il_1(t_2+x_1)}  \phi_{r,s}^{0,0} (-t_1,-t_2) dt_1
 dt_2 \\ &
 =e^{-\frac{\pi(r^2+s^2)}{2q}+\frac{\pi irs}{q}}
 \sum\limits_{l_1,l_2 \operatorname{even}} \ \int\limits_{[0,1]^2}
e^{2\pi i(l_1t_1+l_2t_2)}  \phi_{r,s}^{0,0} (x_2-t_2,x_1-t_1) dt_1  dt_2 \\ &
=\frac{1}{4}\, e^{-\frac{\pi(r^2+s^2)}{2q}+\frac{\pi irs}{q}}
\bigg(  \phi_{r,s}^{0,0} (x_2,x_1) +\phi_{r,s}^{0,0} \bigg( x_2,x_1 -\frac{1}{2}
\bigg) +\phi_{r,s}^{0,0} \bigg( x_2-\frac{1}{2} , x_1 \bigg)  \\ & \qquad \qquad \qquad \qquad \qquad \qquad \qquad \qquad
+\phi_{r,s}^{0,0} \bigg( x_2-\frac{1}{2}\, ,x_1-\frac{1}{2} \bigg) \bigg)  \, .
\end{split}
\end{equation}
Similar computations yield
\begin{equation}
\begin{split}
 S_{01} =\frac{1}{4}\, e^{-\frac{\pi(r^2+s^2)}{2q}+\frac{\pi irs}{q}}
\bigg(  \phi_{r,s}^{0,1} (x_2,x_1) & +\phi_{r,s}^{0,1} \bigg( x_2,x_1 -\frac{1}{2}
\bigg) -\phi_{r,s}^{0,1} \bigg( x_2-\frac{1}{2} , x_1 \bigg) \\ & -
\phi_{r,s}^{0,1} \bigg( x_2-\frac{1}{2}\, ,x_1-\frac{1}{2} \bigg) \bigg)   .
\end{split}
\end{equation}
\begin{equation}
\begin{split}
 S_{10} =\frac{1}{4}\, e^{-\frac{\pi(r^2+s^2)}{2q}+\frac{\pi irs}{q}}
\bigg(  \phi_{r,s}^{1,0} (x_2,x_1) & -\phi_{r,s}^{1,0} \bigg( x_2,x_1 -\frac{1}{2}
\bigg) +\phi_{r,s}^{1,0} \bigg( x_2-\frac{1}{2} , x_1 \bigg) \\  & -
\phi_{r,s}^{1,0} \bigg( x_2-\frac{1}{2} ,x_1-\frac{1}{2} \bigg) \bigg)  .
\end{split}
\end{equation}
\begin{equation}
\begin{split}
 S_{11} =\frac{1}{4}\, e^{-\frac{\pi(r^2+s^2)}{2q}+\frac{\pi irs}{q}}
\bigg( \phi_{r,s}^{1,1} (x_2,x_1)&  - \phi_{r,s}^{1,1} \bigg( x_2,x_1 -\frac{1}{2}
\bigg) -\phi_{r,s}^{1,1} \bigg( x_2-\frac{1}{2} , x_1 \bigg) \\ &
+\phi_{r,s}^{1,1} \bigg( x_2-\frac{1}{2} ,x_1-\frac{1}{2} \bigg) \bigg)   .
\end{split}
\end{equation}
Immediate computations based on the obvious equalities
$$ \vartheta (z,\tau)  =\vartheta^{\operatorname{even}} (z,\tau)+
\vartheta^{\operatorname{odd}} (z,\tau) \qquad \mbox{\rm and} \qquad
 \vartheta_{0,\frac{1}{2}} (z,\tau)  =\vartheta^{\operatorname{even}} (z,\tau) -
\vartheta^{\operatorname{odd}} (z,\tau) $$
show that $\phi_{r,s}^{0,0} =\phi_{r,s}^{0,1} =\phi_{r,s}^{1,0}
=-\phi_{r,s}^{1,1} =\phi_{r,s}$. Actually, for all $r,s\in \Z$, $t_1,t_2 \in\R$,
\begin{equation}
\phi_{r,s} (t_1,t_2) =\frac{\vartheta \big( t_1+\frac{r+is}{2} \big)
\vartheta \big( t_2+\frac{s+ir}{2} \big) -2
\vartheta^{\operatorname{odd}} \big( t_1 +\frac{r+is}{2} \big) \vartheta^{\operatorname{odd}} \big(
t_2 +\frac{s+ir}{2} \big)}{\vartheta (t_1)\, \vartheta (t_2) -2
\vartheta^{\operatorname{odd}} (t_1)  \vartheta^{\operatorname{odd}} (t_2)} .
\end{equation}

From (3.22)-(3.25),
\begin{equation}
\begin{array}{l}
 \displaystyle \sum\limits_{l_1,l_2} \alpha_{ql_1+r,ql_2+s} \,
e^{2\pi i(l_1x_1+l_2x_2)} =S_{00}+S_{01}+S_{10}+S_{11} \\ \\  \displaystyle
 \qquad \qquad =\frac{1}{2}\, e^{-\frac{\pi(r^2+s^2)}{2q}+\frac{\pi irs}{q}} \,
\sum\limits_{\delta_1,\delta_2 =0}^1 \! \! \! (-1)^{\delta_1\delta_2} \,
\phi_{r,s} \bigg( x_2 -\frac{\delta_2}{2} \, ,x_1 -\frac{\delta_1}{2} \bigg) \, .
\end{array}
\end{equation}

We identify $U_1^q$ with $e^{2\pi ix_1}$ and $U_2^q$ with $e^{2\pi ix_2}$
to obtain
\begin{equation}
\begin{split}
 \sum\limits_{l_1,l_2} \alpha_{ql_1+r,ql_2+s} &
U^{ql_2}_2 U^{ql_1}_1
 =\frac{1}{2}\, e^{-\frac{\pi(r^2+s^2)}{2q}+\frac{\pi irs}{q}} \\
& \qquad \times \left( \vartheta \bigg( U_2^q ,\frac{iq}{2} \bigg)
\vartheta \bigg( U_1^q ,\frac{iq}{2} \bigg) -2
\vartheta^{\operatorname{odd}} \bigg( U_2^q ,\frac{iq}{2} \bigg)\
\vartheta^{\operatorname{odd}} \bigg( U_1^q ,\frac{iq}{2} \bigg) \right)^{-1} \\ &
\qquad \times
\sum\limits_{\delta_1,\delta_2 \in \{ 0,1\}} (-1)^{\delta_1\delta_2}
\bigg(  \vartheta \bigg( e^{\pi i(r+\delta_2 +is)} \, U_2^q ,\frac{iq}{2} \bigg)
\vartheta \bigg( e^{\pi i(s+\delta_1 +ir)}  U_1^q ,\frac{iq}{2} \bigg)\\ & \qquad \qquad -
2 \vartheta^{\operatorname{odd}} \bigg( e^{\pi i(r+\delta_2 +is)}  U_2^q ,\frac{iq}{2} \bigg)
\vartheta^{\operatorname{odd}} \bigg( e^{\pi i(s+\delta_1 +ir)}  U_1^q ,\frac{iq}{2} \bigg) \bigg),
\end{split}
\end{equation}
where we denote for any unitary $U$ and $\tau \in \HH$, $q\in {\mathbf N}^*$,
$a\in \Z$, $b\in \R$,
\begin{equation*}
\begin{split}
\vartheta^{\operatorname{odd}} (z,\tau )   & =\sum\limits_{m \operatorname{odd}} e^{\pi i\tau m^2} U^m
=\vartheta^{(2)}_{\frac{1}{2},0}  (U,4\tau)  , \\
\vartheta^{\operatorname{even}} (z,\tau ) & =\sum\limits_{m \operatorname{even}}
e^{\pi i\tau m^2} U^m =\vartheta^{(2)}_{0,0} (U,4\tau)
\end{split}
\end{equation*}
and
\begin{equation*}
\begin{split}
\vartheta^{(q)\operatorname{odd}}_{\frac{a}{q},b} (U,\tau ) & =\sum\limits_{m \operatorname{odd}}
e^{\pi i\tau ( m+\frac{a}{q} )^2 +2\pi i( m+\frac{a}{q} ) b}
U^{qm+a} =\vartheta^{(2q)}_{\frac{a}{2q}+\frac{1}{2},2b} (U,4\tau )  , \\
\vartheta^{(q)\operatorname{even}}_{\frac{a}{q},b} (U,\tau ) & =\sum\limits_{m\operatorname{even}}
e^{\pi i\tau ( m+\frac{a}{q} )^2 +2\pi i( m+\frac{a}{q} ) b}
U^{qm+a} =\vartheta^{(2q)}_{\frac{a}{2q},2b} (U,4\tau ) .
\end{split}
\end{equation*}
As $\vartheta (U,\tau )=\vartheta^{\operatorname{odd}} (U,\tau)+\vartheta^{\operatorname{even}} (U,\tau)$, we get
\begin{equation}
\begin{split}
\vartheta \bigg( U_2^q ,\frac{iq}{2} \bigg) \vartheta \bigg( U_1^q ,\frac{iq}{2}
\bigg)  -2 & \vartheta^{\operatorname{odd}} \bigg( U_2^q ,\frac{iq}{2} \bigg) \vartheta^{\operatorname{odd}}
\bigg( U_1^q ,\frac{iq}{2} \bigg)  \\  &
 =\sum\limits_{\varepsilon_1,\varepsilon_2 =0}^1
(-1)^{\varepsilon_1 \varepsilon_2}
\vartheta^{(2)}_{\frac{\varepsilon_2}{2},0} (U_2^q ,2iq)
\vartheta^{(2)}_{\frac{\varepsilon_1}{2},0} (U_1^q, 2iq) .
\end{split}
\end{equation}
For any unitary $U$ and $n,s\in\Z$, $0\leq s<q$,
\begin{equation*}
\begin{split}
& e^{-\frac{\pi s^2}{2q}} U^s  \vartheta^{\operatorname{even}} \bigg( e^{\pi i(n+is)} U^q ,\frac{iq}{2}\bigg) =
\vartheta^{(q)\operatorname{even}}_{\frac{s}{q},0} \bigg( U,\frac{iq}{2}\bigg)=
\vartheta^{(2q)}_{\frac{s}{2q},0} (U,2iq) ,  \\
& e^{-\frac{\pi s^2}{2q}} U^s \vartheta^{\operatorname{odd}} \bigg( e^{\pi i(n+is)} U^q ,\frac{iq}{2}\bigg) =
e^{\pi in}\vartheta^{(q)\operatorname{odd}}_{\frac{s}{q},0} \bigg( U,\frac{iq}{2}\bigg)= e^{\pi in}
\vartheta^{(2q)}_{\frac{s+q}{2q},0} (U,2iq) ,
\end{split}
\end{equation*}
showing that
\begin{equation}
\begin{split}
e^{-\frac{\pi(r^2+s^2)}{2q}} U_2^s & \bigg( \vartheta \bigg( e^{\pi i(r+\delta_2+is)} U_2^q ,\frac{iq}{2}\bigg) \vartheta
\bigg( e^{\pi i(s+\delta_1+ir)} U_1^q ,\frac{iq}{2} \bigg) \bigg. \\
& \qquad -2\vartheta^{\operatorname{odd}} \bigg( e^{\pi i(r+\delta_2 +is)} U_2^q, \frac{iq}{2}\bigg)
\vartheta^{\operatorname{odd}} \bigg( e^{\pi i(s+\delta_1+ir)} U_1^q ,\frac{iq}{2}\bigg) \bigg) U_1^r     \\
& =\sum\limits_{\varepsilon_1,\varepsilon_2=0}^1 e^{\pi i\varepsilon_1 \varepsilon_2 +\pi i\varepsilon_2 (r+\delta_2)+\pi i\varepsilon_1 (s+\delta_1)}
\vartheta^{(2q)}_{\frac{s+q\varepsilon_2}{2q},0} (U_2,2iq)
\vartheta^{(2q)}_{\frac{r+q\varepsilon_1}{2q},0} (U_1,2iq) .
\end{split}
\end{equation}

We make use of (3.11), (3.28), (3.29), (3.30) and $\sum\limits_{\delta_1,\delta_2 =0}^1
e^{\pi i(\delta_1\delta_2+\varepsilon_1 \varepsilon_2 +\varepsilon_2 \delta_2 +\varepsilon_1\delta_1)} =2$,
$\varepsilon_1,\varepsilon_2=0,1$ to get
\begin{equation*}
e=\frac{\displaystyle \sum\limits_{r,s=0}^{q-1} \sum\limits_{\varepsilon_1,\varepsilon_2=0}^1
e^{\frac{\pi i(r+q\varepsilon_1)(s+q\varepsilon_2)}{q}+\pi i\varepsilon_1\varepsilon_2}
\vartheta^{(2q)}_{\frac{s+q\varepsilon_2}{2q},0} (U_2,2iq)
\vartheta^{(2)}_{\frac{r+q\varepsilon_1}{2q},0} (U_1,2iq)}{\displaystyle
\sum\limits_{\varepsilon_1,\varepsilon_2=0}^1 e^{\pi i\varepsilon_1 \varepsilon_2} \vartheta^{(2)}_{\frac{\varepsilon_2}{2},0}
(U_2^q,2iq) \vartheta^{(2)}_{\frac{\varepsilon_1}{2},0} (U_1^q,2iq)}  ,
\end{equation*}
which implies the formula in Proposition 3.1 (ii).
\end{proof}

\medskip

We continue with some considerations on the case
 $\alpha = 2^{-1} .$ According to Proposition 3.1, the formula
$$e=\frac{1}{2}
\left( 1+\frac{\vartheta^{(2)}_{0,\frac{1}{2}} (U_2,i)\,
 \vartheta^{(2)}_{\frac{1}{2},0} (U_1,i)}{
\vartheta (U_2^2,i) \, \vartheta (U_1^2,i)}+
\frac{\vartheta^{(2)}_{\frac{1}{2},0} (U_2,i) \,
 \vartheta^{(2)}_{0,\frac{1}{2}} (U_1,i)}{
\vartheta (U_2^2,i)\, \vartheta (U_1^2,i)} -\frac{
i\, \vartheta^{(2)}_{\frac{1}{2},\frac{1}{2}} (U_2,i)\,
\vartheta^{(2)}_{\frac{1}{2},
\frac{1}{2}} (U_1,i)}{\vartheta (U_2^2,i)\, \vartheta (U_1^2,i)}
\right)
$$
defines a
projection of trace $\frac{1}{2}$
 in $A_{\frac{1}{2}}$. We shall simply
denote $e =e_1 =( + \ + \ + \ -)\, .$

We consider the automorphisms
$\rho_{t_1,t_2}$ of $A_{\frac{1}{2}}$ acting on the
generators $U_1$ and $U_2$ by $\rho_{t_1,t_2} (U_j)=e^{2\pi it_j} U_j $, $j=1,2$,
$t_1 ,t_2 \in \R$ and $\rho$ acting by $\rho (U_1)=-U_2$, $\rho (U_2)=-U_1$.
 The relations
\begin{equation*}
\begin{split}
&  \vartheta^{(2)}_{\frac{1}{2},b} ( -U,\tau ) =
-\vartheta^{(2)}_{\frac{1}{2},b} (U,\tau) ,\\  &
 \vartheta^{(2)}_{\frac{1}{2},\frac{1}{2}} (U_1,\tau )
\vartheta^{(2)}_{\frac{1}{2},\frac{1}{2}} (U_2,\tau )=-
\vartheta^{(2)}_{\frac{1}{2},\frac{1}{2}} (U_2,\tau )
\vartheta^{(2)}_{\frac{1}{2},\frac{1}{2}} (U_1,\tau ),\qquad b\in\R,\ \tau\in \HH ,
\end{split}
\end{equation*}
and the fact that $\vartheta (U_1^2,i)$ and
$\vartheta (U_2^2 ,i)$ are central in $A_{\frac{1}{2}}$ show that
\begin{equation*}
\begin{split}
& e_1  =(+\ +\ +\ -) , \qquad \quad
e_2 =\rho_{\frac{1}{2},0} (e_1) =(+\ -\ +\ +) , \qquad \quad
e_3 =\rho_{0,\frac{1}{2}} (e_1) =(+\ +\ -\ +)  , \\ &
e_4 =\rho_{\frac{1}{2},\frac{1}{2}} (e_1) =(+\ -\ -\ -)  , \qquad \quad
e_5 =\rho(e_1) =(+\ -\ -\ +)  .
\end{split}
\end{equation*}
In particular $e_j$, $1\leq j\leq 5$ are projections of trace $\frac{1}{2}$
in $A_{\frac{1}{2}}$ such that
\begin{equation}
 e_1 +e_5 =e+\rho (e) =1  ,
\end{equation}
\begin{equation}
 e_1 +e_2 +e_3 +e_4 =\sum\limits_{c_1,c_2=0}^1 \rho_{\frac{c_1}{2},\frac{c_2}{2}} (e)=2 .
\end{equation}
We also notice that if $\tilde{\rho}$
 is the automorphism of $A_{\frac{1}{2}}$ defined by $\tilde{\rho} (U_1)=U_2$,
$\tilde{\rho} (U_2)=U_1$, then
$$\tilde{\rho} (e_1 )=(+\ +\ +\ +).$$

Actually, the analogue of (3.32) holds
for any $\alpha =q^{-1}$, $q\in {\mathbf N}^*$, $q\geq 2$. To see this,
we notice first
that for all $s,c \in \Z$, $b\in \R$,
$$\vartheta^{(q)}_{\frac{s}{q},b} \big( e^{\frac{2\pi ic}{q}} U ,\tau \big) =
e^{\frac{2\pi isc}{q}}\vartheta^{(q)}_{\frac{s}{q},b} (U ,\tau ) ,$$
so if we write $e=\sum\limits_{r,s=0}^{q-1} e^{-\frac{\pi irs}{q}} A_{r,s}$ as in Proposition 3.1, then
$$\rho_{\frac{c_1}{q},\frac{c_2}{q}} (e)=\frac{1}{q} \sum\limits_{r,s=0}^{q-1}
e^{-\frac{\pi irs}{q}+\frac{2\pi isc_2}{q} +\frac{2\pi irc_1}{q}}  A_{r,s} ,$$
and furthermore
\begin{equation}
\sum\limits_{c_1,c_2 =0}^{q-1} \rho_{\frac{c_1}{q},\frac{c_2}{q}} (e)
=\frac{1}{q} \sum\limits_{r,s=0}^{q-1} e^{-\frac{\pi irs}{q}}
\sum\limits_{c_1 =0}^{q-1} e^{\frac{2\pi irc_1}{q}}
\sum\limits_{c_2 =0}^{q-1} e^{\frac{2\pi isc_2}{q}}   A_{r,s}
=qA_{0,0} =q  .
\end{equation}

Let $\pi :A_{\frac{1}{q}} \rightarrow M_q (\C)$, $\pi (U_j )=\tilde{U}_j$,
$j=1,2$ be the canonical finite dimensional representation of
 $A_{\frac{1}{q}}$ (i.e. $\tilde{U}_1^q =\tilde{U}_2^q =I_q$). We
consider the finite dimensional representations
$\pi_{t_1,t_2} =\pi  \rho_{t_1,t_2}$ of $A_{\frac{1}{q}}$ and get
$\pi_{t_1 ,t_2} (U_j^{qm+s} )=e^{2\pi i(qm+s)t_j}\, \tilde{U}_j^s$, $j=1,2$,
$0\leq s<q$. Furthermore, we obtain
\begin{equation*}
\begin{split}
& \pi_{t_1,t_2} \left(   \vartheta^{(q)}_{\frac{s}{q},\frac{r}{2}}
(U_2 ,\tau ) \right) =\sum\limits_m e^{\pi i\tau ( m+\frac{s}{q} )^2
+2\pi i ( m+\frac{s}{q} ) \frac{r}{2}}
e^{2\pi i(qm+s)t_2} \tilde{U}_2^s =
\vartheta_{\frac{s}{q},\frac{r}{2}} (qt_2 ,\tau )\, \tilde{U}^s_2  ,
\\  &
\pi_{t_1,t_2} \left( \vartheta^{(q)}_{\frac{r}{q},\frac{s}{2}}
(U_1 ,\tau ) \right)
 =\vartheta_{\frac{r}{q},\frac{s}{2}} (qt_1 ,\tau )\, \tilde{U}^r_1  ,\\
&
\pi_{t_1,t_2} \big( \vartheta (U_j^q ,\tau )\big)  =\vartheta (qt_j ,\tau )
\cdot I_q ,\quad j=1,2 .
\end{split}
\end{equation*}

Therefore, for any non-negative even integer $q$,
\begin{equation}
\pi_{t_1,t_2} (e)=q^{-1} \vartheta \bigg( qt_2 ,\frac{iq}{2} \bigg)^{-1}
\vartheta \bigg( qt_1 , \frac{iq}{2} \bigg)^{-1}
\sum\limits_{r,s=0}^{q-1} e^{-\frac{\pi irs}{q}}
\vartheta_{\frac{s}{q},\frac{r}{2}} \bigg( qt_2 ,\frac{iq}{2} \bigg)
\vartheta_{\frac{r}{q},\frac{s}{2}} \bigg( qt_1 ,\frac{iq}{2} \bigg)
\tilde{U}_2^s \tilde{U}_1^r.
\end{equation}

We compute $\tau \big( \pi_{t_1,t_2} (e)^2 \big)$ using (3.34) and obtain
$\tau \big( \pi_{t_1,t_2} (e)\big) =q^{-1}$. Therefore, for all $t_1,t_2 \in \R$,
\begin{equation}
\begin{split}
 \sum\limits_{m,n=0}^{q-1}
 e^{-\frac{4\pi imn}{q}}   \vartheta_{-\frac{n}{q},-\frac{m}{2}} &
\bigg( t_1 ,\frac{iq}{2} \bigg)
\vartheta_{\frac{n}{q},\frac{m}{2}}
\bigg( t_1 ,\frac{iq}{2} \bigg)
\vartheta_{-\frac{m}{q}, -\frac{n}{2}} \bigg( t_2 ,\frac{iq}{2} \bigg)
\vartheta_{\frac{m}{q}, \frac{n}{2}} \bigg( t_2 ,\frac{iq}{2} \bigg) \\ &
=q \vartheta \bigg( t_1 ,\frac{iq}{2} \bigg)^2
\bigg( t_2 ,\frac{iq}{2} \bigg)^2  .
\end{split}
\end{equation}

Using $\vartheta_{-a,-b} (z,\tau)=\vartheta_{a,b} (-z,\tau)$,
it follows that for any even integer $q\in {\mathbf N}^*$ and any $t_1,t_2 \in
\R$,
\begin{equation}
\begin{split}
 \sum\limits_{m,n=0}^{q-1}  e^{-\frac{4\pi imn}{q}} &
 \vartheta_{\frac{n}{q},\frac{m}{2}} \bigg( t_1 ,
\frac{iq}{2} \bigg)
 \vartheta_{\frac{n}{q},\frac{m}{2}} \bigg(- t_1 ,
\frac{iq}{2} \bigg)
\vartheta_{\frac{m}{q},\frac{n}{2}} \bigg( t_2 ,\frac{iq}{2} \bigg)
\vartheta_{\frac{m}{q},\frac{n}{2}} \bigg( - t_2 ,\frac{iq}{2} \bigg)  \\ &
 =q \vartheta \bigg( t_1 ,\frac{iq}{2} \bigg)^2 \vartheta\bigg( t_2 ,\frac{iq}{2} \bigg)^2  .
\end{split}
\end{equation}

Actually (3.38) can be regarded as a sort of Riemann theta relation.
To see this, we specialize further to $q=2$ and denote as
in \cite{M} $\vartheta_{00}=\vartheta$, $\vartheta_{10} =\vartheta_{\frac{1}{2},0}$,
$\vartheta_{01} =\vartheta_{0,\frac{1}{2}}$, $\vartheta_{11} =\vartheta_{\frac{1}{2},
\frac{1}{2}}$. Taking $q=2$ in (3.38) and using the obvious equalities
\begin{equation*}
\begin{split}
& \vartheta_{00} (-z,\tau) =\vartheta_{00} (z,\tau) ,
\qquad
\vartheta_{10} (-z,\tau) =\vartheta_{10} (z,\tau) , \\ &
 \vartheta_{01} (-z,\tau) =\vartheta_{01} (z,\tau) ,\qquad
\vartheta_{11} (-z,\tau )=-\vartheta_{11} (z,\tau)  ,
\end{split}
\end{equation*}
we obtain for all $x,u\in \R$,
\begin{equation}
 \vartheta_{01} (x,i)^2  \vartheta_{10} (u,i)^2 +
\vartheta_{10} (x,i)^2  \vartheta_{01} (u,i)^2 +
\vartheta_{11} (x,i)^2  \vartheta_{11} (u,i)^2 =
\vartheta_{00} (x,i)^2 \vartheta_{00} (u,i)^2 .
\end{equation}

\begin{remark}
{\em Identity (3.37) should be compared with
 Riemann's theta formulae, for example with
formula (A$_1$) in \cite[p.21]{M}}
\begin{equation*}
\begin{split}
\vartheta_{00} (x,\tau)^2  \vartheta_{00} (u,\tau )^2 +
\vartheta_{11} (x,\tau )^2 \vartheta_{11} (u,\tau)^2
 & =
\vartheta_{01} (x,\tau)^2 \vartheta_{01} (u,\tau )^2 +
\vartheta_{10} (x,\tau )^2 \vartheta_{10} (u,\tau )^2 \qquad \quad \mbox{\em (A$_1$)}  \\ &
 =\vartheta_{00} (x+u,\tau ) \vartheta_{00} (x-u,\tau ) \vartheta_{00} (0,\tau )^2 .
\end{split}
\end{equation*}
\end{remark}

\appendix
\numberwithin{equation}{section}
\makeatletter
\makeatother

\section{}

Let $\alpha \in [0,1)$, $\rho =e^{2\pi i\alpha}$, $A_\alpha =C^* (u,v; uv=\rho vu)$ and $\sigma$ be the order
four automorphism of $A_\alpha$ defined by $\sigma (u)=v$, $\sigma (v)=
u^{-1}$. Denote by $E$ the conditional expectation from $A_\alpha$ onto
$A_\alpha^\sigma$ defined by
$$E(x)=\frac{1}{4}   \big( x+\sigma (x)+\sigma^2 (x) +
\sigma^3 (x)\big),   \qquad x\in A_\alpha   .$$

For any $n,m\in \Z$, set
\begin{equation*}
\begin{split}
[ n,m ] &  =\rho^{-\frac{nm}{2}} (u^n v^m +u^{-n} v^{-m} )
=2\, \rho^{-\frac{nm}{2}} \big( u^n v^m +\sigma^2 (u^n v^m )\big) ,\\
 \{ n,m\} & = 4 \rho^{-\frac{nm}{2}} E(u^n v^m )=
\rho^{-\frac{nm}{2}} (u^n v^m +u^{-n} v^{-m} +v^n u^{-m}
+v^{-n} u^m )  \\  &
 =\rho^{-\frac{nm}{2}} (u^n v^m +u^{-n} v^{-m} )+
\rho^{\frac{nm}{2}} (u^{-m} v^n +u^m v^{-n} )
 =[  n,m]  +[ -m ,n]   .
\end{split}
\end{equation*}

The following properties of $[ n,m ] $ are easy to check (\cite{BEEK}):
\begin{equation}
 [n,m]^* = [n,m]  =[-n,-m]   ,
\end{equation}
\begin{equation}
[n,m] [k,l]   =\rho^{\frac{nl-mk}{2}} [n+k,m+l]  +
\rho^{\frac{mk-nl}{2}} [n-k,m-l]  .
\end{equation}

Moreover, using (A.1) and (A.2) it is plain to check that for all $n,m,k,l\in \Z$,
\begin{equation}
 \{ n,m\}^\ast = \{ n,m\} =\{ -m,n\}  = \{ -n,-m\} =\{ m,-n\} ,
\end{equation}
\begin{equation}
\begin{split}
\{ n,m\}  \{ k,l\} &  =  \rho^{\frac{nl-mk}{2}} \{
n+k,m+l \} +\rho^{\frac{mk-nl}{2}} \{ n-k,m-l \}  \\ &  +
\rho^{\frac{nk+ml}{2}} \{ n-l ,m+k \}   +
\rho^{\frac{-nk-ml}{2}} \{ n+l ,m-k \} .
\end{split}
\end{equation}

\begin{prop}
If $\alpha \neq 0$ and $\alpha \neq \frac{1}{2}$,  then
$A_\alpha^\sigma =C^* \big( \{ n,0\}  ;\ n\geq 0 \big) .$
\end{prop}

\begin{proof}
The linear span of $(u^n v^m )_{n,m\in \Z}$ is dense in $A_\alpha$ and
the conditional expectation
$E$ continuous, hence $\operatorname{span} \big( \{ n,m\} \big)_{n,m \in \Z}$ is dense
in $A_\alpha^\sigma$. Denote by $B$
the $C^*$-algebra generated by $\big( \{ n,0\}
\big)_{n\in {\mathbf N}}$. According to (A.3), it suffices to prove
that for all $n,m\in {\mathbf N}$,
\begin{equation}
\{ n,m\} \in B.
\end{equation}
According to (A.4) and (A.3)
\begin{equation}
\{ n,0 \} \{ k,0\} =\{ n+k ,0\} +\{ n-k,0\} +\rho^{\frac{nk}{2}} \{ n,k\}
+\rho^{-\frac{nk}{2}} \{ k,n\}  ,
\end{equation}
whence $\{ n,k\} \in B$ if and only if $\{ k,n\} \in B$.
The previous equality gives also for all $n,k\in \mbox{\bf Z}$
\begin{equation}
\rho^{\frac{nk}{2}} \{ n,k\} +
\rho^{-\frac{nk}{2}} \{ k,n\} \in B .
\end{equation}

In particular $\{ 1,1\}  \in B$ and using again (A.4) and (A.3) we get for all $n,k\in\Z$,
\begin{equation}
 \{ n,0\} \{ 1,1\} =\rho^{\frac{n}{2}}
\{ n+1,1\} +\rho^{-\frac{n}{2}} \{ 1,n-1\} +
\rho^{\frac{n}{2}} \{ n-1,1\} +
\rho^{-\frac{n}{2}} \{ 1,n+1\}  .
\end{equation}

Let $m\geq 1$ and assume that $\{ k,1\} \in B$ (hence also $\{ 1,k\}\in B$) for all $
0\leq k\leq m$. By (A.8)
\begin{equation}
\rho^{\frac{m}{2}} \{ m+1,1\}  +
\rho^{-\frac{m}{2}} \{ 1,m+1\} \in B .
\end{equation}

Taking $k=1$, $n=m+1$ in (A.7) and using (A.9) and
$$\left| \begin{matrix}  \rho^{\frac{m+1}{2}} & \rho^{-\frac{m+1}{2}} \\
\rho^{\frac{m}{2}} & \rho^{-\frac{m}{2}} \end{matrix}
\right|  =\rho^{\frac{1}{2}} -
\rho^{-\frac{1}{2}} \neq 0 ,$$
we conclude that $\{ m+1,1\} \in B$, hence $\{ n,1\} \in B$ and
$\{ 1,n\} \in B$ for all $n\geq 0$.

Finally, we prove (A.5) by induction for $n,m\geq 0$. Assume that for some $k\geq 1$ we have
\begin{equation}
\{ n,m\} \in B\ \  \mbox{for all $(n,m)\in \big( [0,\infty) \times [0,k]\big)
\times \big( [0,k]\times [0,\infty)\big)$}.
\end{equation}

To conclude, it will suffice, according to (A.7), to prove that for all $n\geq 0$,
$$\{ n,k+1 \} \in B  .$$

This holds for $n=0,1,\ldots,k$. By (A.10) and
\begin{equation*}
\begin{split}
\{ k,k\} \{ 1,1\}  &
=\{ k+1,k+1 \} +\{ k-1,k-1\} +\rho^k \{ k-1,k+1 \} +\rho^{-k} \{ k+1,k-1\} \in B , \\
\{ k+1,k\}  \{ 1,1\} &
=\rho^{\frac{1}{2}} \{ k+2,k+1 \}  +
\rho^{-\frac{1}{2}} \{ k,k-1\}  + \rho^{k+\frac{1}{2}} \{ k,k+1\}
 + \rho^{-k-\frac{1}{2}} \{ k+2,k-1\} \in B ,  \\
\{ k+2,k\} \{ 1,1\} &  =\rho \{ k+3,k+1\}  +
\rho^{-1} \{ k+1 ,k-1\}  +
\rho^{k+1} \{ k+1,k+1\} \\ & \quad
 + \rho^{-k-1} \{ k+3,k-1\}  \in  B, \quad \mbox{\rm etc.},
\end{split}
\end{equation*}
we get $\{ k+1,k+1\} ,\{ k+2,k+1\} ,\{ k+3,k+1\} ,\ldots \in B$.
\end{proof}

\begin{cor}
If $\alpha \notin \mbox{\bf Q}$, then
$$A_\alpha^\sigma =C^* \big( \{ 1,0\} , \{ 2,0\} \big) =
C^* \big( u+v+u^{-1} +v^{-1} ,u^2 +v^2 +u^{-2} +v^{-2} \big) \, .$$
\end{cor}

\begin{proof}
Denote by $C$ the $\ast$-algebra generated by $\{ 1,0\} $ and
$\{ 2,0 \}$. According to the previous proposition it suffices to
show that $\{ n,0\} \in C$ for all $n\geq 3$. Firstly,
notice that
$$\{ 1,0\}^2 = \{ 2,0\}  +\{ 0,0\} +
\big( \rho^{\frac{1}{2}} +\rho^{-\frac{1}{2}} \big) \,
 \{ 1,1\} \, \in \,  C$$
yields $\{ 1,1\} \in C$. Then, prove by induction that for all $n\geq 2$,
\begin{equation}
\{ 1,0\} ,\ldots ,\{ n,0 \} ,\{ 1,1\} ,\ldots , \{ n-1 ,1\} \in C .
\end{equation}

Assume that (A.11) holds for some $n\geq 2$. By (A.6),
$$\rho^{\frac{n-2}{2}} \{ n-2,1\}  +\rho^{-\frac{n-2}{2}}  \{ 1,n-2\}
 = \{ n-2,0\}  \{ 1,0\}  - \{ n-1,0\}  - \{ n-3,0\}  \in C ,$$
hence $ \{ 1,n-2\} =\{ 2-n,1\} \in C$. By (A.8) and the induction hypotheses
\begin{equation}
\rho^{\frac{n-1}{2}} \{ n,1\} +\rho^{-\frac{n-1}{2}}
\{ 1,n\} =\{ n-1,0\}  \{ 1,1\}
-\rho^{-\frac{n-1}{2}} \{ 1,n-2\} -\rho^{\frac{n-1}{2}} \{ n-2,1 \}\in C  .
\end{equation}

On the other hand
$$ \{ 1,1\}  \{ n-1,0 \} =\rho^{-\frac{n-1}{2}}
 \{ n,1 \} +\rho^{\frac{n-1}{2}}  \{ 1,n\}  +
\rho^{\frac{n-1}{2}} \{ 2-n,1\}  +
\rho^{-\frac{n-1}{2}}  \{ n-2,1 \} \in C  ,$$
whence
\begin{equation}
\rho^{-\frac{n-1}{2}} \{ n,1\} +
\rho^{\frac{n-1}{2}}  \{ 1,n\} \in C  .
\end{equation}

Since $\rho^m \neq 1$, $m\in\Z$, (A.11) and (A.12) show that
$\{ n,1\} ,  \{ 1,n \} \in C$. Furthermore, using also
$$\{ n,0 \} \{ 1,0\}  =\{ n+1,0\} +
\{ n-1,0 \} +\rho^{\frac{n}{2}} \{ n,1\} +
\rho^{-\frac{n}{2}} \{ 1,n \} \in C ,$$
we conclude that $ \{ n+1,0\} \in C$. Summarizing, we have shown that $\{ n,0\}, \{ n,1\} \in C$ for all $n\geq 0$.

Computing $\{ 1,0\} \{ 1,1\}, \{ 1,1\}^2, \{ 1,2\} \{ 1,1\} ,\ldots $ we show that
$\{ n,2\} \in C$ for all $n\geq 0$. Finally, one checks by induction on $k\geq 0$ as above that
$\{ n,k\} \in C$ for all $n\geq 0$.
\end{proof}

\bigskip

{\sl Addendum.} {\small After writing this paper, I have received the preprint
``Chern characters of Fourier modules" by S. Walters, where the author computes the Chern characters
of nine basic modules over the crossed product $A_\alpha \propto_\sigma \Z_4$. Those computations involve also Jacobi's theta functions.
I was recently informed by A. Valette that estimates on the norm of the Harper operator have been obtained by different
methods in the papers ``On the spectrum of a random walk on the discrete Heisenberg group and the norm of Harper's operators",
J. Geom. Phys. {\bf 21}, 337-356 (1997), by C. B\' eguin, A. Valette, A. Zuk and in
``Norm estimates of discrete Schr\" odinger operators", Colloq. Math. {\bf 76}, 153-160 (1998) by R. Szwarc.}

\bigskip

{\sl Acknowledgment.} {\small I would like to express my gratitude to Dai Evans for his continuous encouragement and
support. I am grateful to Bernard Helffer for explaining me the results from \cite{HS}.}

\end{document}